\documentclass[oneside, letterpaper]{amsart}

\usepackage{amsmath, amssymb, bm}
\usepackage{enumitem}
\usepackage{xy}
\xyoption{all}
\usepackage[letterpaper]{geometry}

\usepackage{setspace}
\usepackage{cite}
\usepackage{url}
\usepackage[usenames, dvipsnames]{xcolor}

\usepackage{amsthm}

\newtheorem{tm}{Theorem}[section]
\newtheorem{theorem}[tm]{Theorem}
\newtheorem{pr}[tm]{Proposition}

\theoremstyle{plain} 
\newcommand{\thistheoremname}{}
\newtheorem{genericthm}[tm]{\thistheoremname}
\newenvironment{namedthm}[1]
  {\renewcommand{\thistheoremname}{#1}%
   \begin{genericthm}}
  {\end{genericthm}}

\theoremstyle{definition}
\newtheorem{df}[tm]{Definition}
\newtheorem{definition}[tm]{Definition}

\theoremstyle{remark}

\newtheorem{example}[tm]{Example}

\newcommand{\sh}[1]{{\mathcal{#1}}}
\newcommand{\cat}[1]{{\mathbf{#1}}}

\newcommand{\isom}{\cong}
\newcommand{\homeo}{\approx}
\newcommand{\iso}{\isom}

\newcommand{\weq}{\simeq}

\newcommand\rightthreearrow{%
        \mathrel{\vcenter{\mathsurround0pt
                \ialign{##\crcr
                        \noalign{\nointerlineskip}$\rightarrow$\crcr
                        \noalign{\nointerlineskip}$\rightarrow$\crcr
                        \noalign{\nointerlineskip}$\rightarrow$\crcr
                }%
        }}%
}

\newcommand{\A}{{\bbb{A}}}
\newcommand{\Aone}{{\mathbb{A}^{\!1}}}

\newcommand{\KMW}{\mathbf{K}^{\mathrm{MW}}}
\newcommand{\KM}{\mathbf{K}^{\mathrm{M}}}
\newcommand{\KMf}{\mathrm{K}^{\mathrm{M}}}
\newcommand{\GW}{\mathrm{GW}}

\newcommand{\Sheafproj}{\underline{\operatorname{Proj}}}

\newcommand{\Nis}{\mathrm{Nis}}

\newcommand{\bpi}{\bm{\pi}}

\DeclareMathOperator*{\hocolim}{hocolim}

\newcommand{\rH}{\operatorname{H}}
\newcommand{\Hoh}{\rH}

\newcommand{\et}{\text{\'et}}
\newcommand{\equi}{\textrm{equi}}
\newcommand{\gm}{\text{gm}}

\newcommand\bbb[1]{\ensuremath{{\mathbb{#1}}}}
\newcommand{\Z}{\mathbb{Z}}

\newcommand{\R}{\mathbb{R}}
\newcommand{\RR}{\R}

\newcommand{\C}{\mathbb{C}}
\newcommand{\G}{\mathbb{G}}
\newcommand{\Gm}{\mathbb{G}_m}
\newcommand{\ZZ}{\Z}

\newcommand{\CC}{\C}
\newcommand{\PP}{\mathbb{P}}

\newcommand{\calO}{\mathcal{O}}

\newcommand{\Map}{\operatorname{Map}} 
\newcommand{\Fun}{\operatorname{Fun}}  
\newcommand{\Mor}{\operatorname{Mor}} 

\newcommand{\Tr}{\operatorname{Tr}}
\newcommand{\Gal}{\operatorname{Gal}}
\newcommand{\Bl}{\operatorname{Bl}}

\newcommand{\Set}{\cat{Set}}
\newcommand{\sSet}{\cat{sSet}}

\newcommand{\Sm}{\cat{Sm}}
\newcommand{\Spaces}[1][{k}]{\cat{sPre}(\Sm_{#1})} 
\newcommand{\Top}{\mathbf{Top}}
\newcommand{\Spt}{\cat{Spt}} 
\newcommand{\soneSpt}{\Spt(\Sm_k)}

\newcommand{\Gr}{\operatorname{Gr}}

\newcommand{\disc}{\operatorname{Disc}}
\newcommand{\ind}{\operatorname{ind}}
\newcommand{\rk}{\operatorname{rank}}
\newcommand{\Jac}{\operatorname{Jac}}
\newcommand{\Th}{\operatorname{Th}}
\newcommand{\Spec}{\operatorname{Spec}}
\newcommand{\GL}{\operatorname{GL}}
\newcommand{\SL}{\operatorname{SL}}
\newcommand{\Sp}{\operatorname{Sp}}
\newcommand{\Ker}{\operatorname{Ker}}

\newcommand{\tensor}{\otimes}
\newcommand{\pt}{\ast}

\newcommand{\op}{{\mathrm{op}}}

\newcommand{\sgn}{\operatorname{sign}}

\newcommand{\Sing}{\operatorname{Sing}}
\newcommand{\Pic}{\operatorname{Pic}}
\newcommand{\Type}{\operatorname{Type}}
\newcommand{\Zar}{\text{Zar}}

\newcommand{\EKL}{\mathrm{EKL}}

\newcommand{\K}{\mathbf{K}}

\newcommand{\hidden}[1]{\footnote{Hidden:  #1}}
\renewcommand{\hidden}[1]{}

\input{figure_header.tex}
\usepackage{graphicx}

\begin{document}
\title{Unstable motivic homotopy theory}
\author{Kirsten Wickelgren}
\address{School of Mathematics, Georgia Institute of Technology, Atlanta~GA, USA}
\email{wickelgren@post.harvard.edu}
\author{Ben Williams}
\address{Department of Mathematics, University of British Columbia, Vancouver~BC, Canada}
\email{tbjw@math.ubc.ca}

\maketitle



\section{Introduction}

Morel--Voevodsky's $\Aone$-homotopy theory transports tools from algebraic topology into arithmetic and algebraic
geometry, allowing us to draw arithmetic conclusions from topological arguments. Comparison results between classical
and $\Aone$-homotopy theories can also be used in the reverse direction, allowing us to infer topological results from
algebraic calculations. For example, see the article by Isaksen and {\O}stv{\ae}r on Motivic Stable Homotopy Groups \cite{Isaksen-Ostvaer-handbook}. The present article will introduce unstable $\Aone$-homotopy theory and give several applications.

Underlying all $\Aone$-homotopy theories is some category of schemes. A special case of a scheme is that of an \textit{affine
 scheme}, $\Spec R$, which is a topological space, the points of which are the prime ideals of a ring $R$ and on which
the there is a sheaf of rings essentially provided by $R$ itself. For example, when $R$ is a finitely generated $k$-algebra, $R$ can be written as $k[x_1, \ldots,x_n]/\langle f_1,\ldots,f_m\rangle$, and $\Spec R$ can be thought of as the common zero locus of the polynomials $f_1$,$f_2$,...,$f_m$, that is to say, $\{ (x_1, \ldots, x_n) : f_i(x_1,\ldots,x_n) = 0 \text{ for } i =1,\ldots,m\}$. Indeed, for a $k$-algebra $S$, the set $$(\Spec R) (S) := \{ (x_1, \ldots, x_n) \in S^n: f_i(x_1,\ldots,x_n) = 0 \text{ for } i =1,\ldots,m\}$$ is the set of $S$-points of $\Spec R$, where an $S$-point is a map $\Spec S \to \Spec R$. We remind the reader that a scheme $X$ is a locally
ringed space that is locally isomorphic to affine schemes. So heuristically, a scheme is formed by gluing together pieces, each of which is the common zero locus of a set of polynomials. The topology on a scheme $X$ is called the \textit{Zariski
 topology}, with basis given by subsets $U$ of affines $\Spec R$ of the form $U = \{ \mathfrak{p} \in \Spec R : g \notin \mathfrak{p} \}$ for some $g$ in $R$. It is too coarse to be of use for classical homotopy theory; for instance, the topological space of an
irreducible scheme is contractible, having the generic point as a deformation retract. Standard references for the theory of schemes
include \cite{hartshorneAlgebraicGeometry1977}, \cite{vakilRisingSeaFoundations2015}, and for the geometric view that
motivates the theory \cite{eisenbudGeometrySchemes2000}.

Both to avoid certain pathologies and to make use of technical theorems that can be proved under certain assumptions,
$\Aone$-homotopy theory restricts itself to considering subcategories of the category of all schemes. In the seminal
\cite{morelMathbbAHomotopyTheory1999} the restriction is already made to $\Sm_S$, the full subcategory of smooth schemes of
finite type over a finite dimensional noetherian base scheme $S$. The smoothness condition, while technical, is
geometrically intuitive. For example, when $S = \Spec \bbb{C}$, the smooth $S$-schemes are precisely those whose $\bbb{C}$-points form a manifold. The ``finite-type'' condition is most easily understood if $S= \Spec R$ is affine, in which
case the finite-type $S$-schemes are those covered by finitely many affine schemes, each of which is determined by the vanishing
of finitely many polynomials in finitely many variables over the ring $R$. The most important case, and the best understood,
is when $S= \Spec k$ is the spectrum of a field. In this case $\Sm_k$ is the category of smooth $k$-varieties.

In this article, we will use the notation $S$ to indicate a base scheme, always assumed noetherian and finite
dimensional. For the more sophisticated results in the sequel, it will be necessary to assume $S = \Spec k$ for some
field. We have not made efforts to be precise about the maximal generality of base scheme $S$ for which a particular
result is known to hold.

To see the applicability of $\Aone$-homotopy theory, consider the following example. The topological Brouwer degree
map \[\deg: [S^n, S^n] \to \Z \] from pointed homotopy classes of maps of the $n$-sphere to itself can be evaluated on a
smooth map \[f: S^n \to S^n\] in the following manner. Choose a regular value $p$ in $S^n$ and consider its finitely
many preimages \hidden{ There are finitely many preimages because the fiber is compact and discrete. See Milnor's {\em
 Topology from the differentiable viewpoint} p8} $f^{-1}(p) = \{ q_1, \ldots, q_m \}$. At each point $q_i$, choose
local coordinates compatible with a fixed orientation on $S^n$. The induced map on tangent spaces can then be viewed as
an $\R$-linear isomorphism $T_{q_i} f: \R^n \to \R^n$, with an associated Jacobian determinant
$\Jac f (q_i) = \det (\frac{\partial (T_{q_i} f)_j}{ x_k})_{j,k}$. The assumption that $p$ is a regular value implies
that $\Jac f (q_i) \neq 0$, and therefore there is a local degree $\deg_{q_i} f$ of $f$ at $q_i$ such that
$$ \deg_{q_i} f = \begin{cases}
 +1 \text{ if } \Jac f (q_i) >0 ,\\
 -1 \text{ if } \Jac f (q_i) < 0.
\end{cases}$$ The degree $\deg f$ of $f$ is then given \[\deg f = \sum_{q \in f^{-1} (p)} \deg_q f,\] as the appropriate sum of $+1$'s and $-1$'s. 

Lannes and Morel suggested the following modification of this formula to give a degree for an algebraic function
$f: \PP^1 \to \PP^1$ valued in nondegenerate symmetric bilinear forms. Namely, let $k$ be a field and let
$\GW(k)$ denote the Grothendieck--Witt group, whose elements are formal differences of $k$-valued, non-degenerate, symmetric, bilinear forms on finite dimensional $k$-vector spaces. We will say more about
this group in Section \ref{degree_section}. For $a$ in $k^\times/(k^\times)^2$, denote by $\langle a \rangle$ the element of
$\GW(k)$ determined by the isomorphism class of the bilinear form $(x,y) \mapsto a xy$. For simplicity, assume that $p$
is a $k$-point of $\bbb{P}^1$, i.e., $p$ is an element of $k$ or $\infty$, and that the points $q$ of $f^{-1}(p)$ are
also $k$-points such that $\Jac f (q) \neq 0$. Then the $\Aone$-degree of $f$ in $\GW(k)$ is given by
\[\deg f = \sum_{q \in f^{-1} (p)} \langle \Jac f (q) \rangle.\] In other words, the $\Aone$-degree counts the points of the inverse
image weighted by their Jacobians, instead of weighting only by the signs of their Jacobians.

Morel shows that this definition extends to an $\Aone$-degree homomorphism \begin{equation}\label{deghom}\deg:
 [\bbb{P}_k^n/\bbb{P}_k^{n-1}, \bbb{P}_k^n/\bbb{P}_k^{n-1}] \to \GW(k)\end{equation} from $\Aone$-homotopy classes of
endomorphisms of the quotient $\bbb{P}_k^n/\bbb{P}_k^{n-1}$ to the Grothendieck--Witt group. We will use this degree to enrich results in classical enumerative geometry over $\bbb{C}$ to equalities in
$\GW(k)$ in section \ref{Applications_to_enumerative_geometry_section}. For \eqref{deghom} to make sense, we must
define the quotient $\bbb{P}_k^n/\bbb{P}_k^{n-1}$. It is not a scheme; it is a space in the sense of
Morel--Voevodsky. In section \ref{Section:construction_unstable_A1-homotopy_theory}, we sketch the construction of the
homotopy theory of spaces, including Thom spaces and the Purity Theorem. In section \ref{RealizationsSection}, we discuss realization functors to topological spaces, which allow us to see how $\bbb{A}^1$-homotopy theory combines phenomena associated to the real and complex points of a variety. An example of this is the degree, discussed in the follow section, along with Euler classes. The Milnor--Witt $K$-theory groups are also introduced in section \ref{degree_section}. These are the global sections of certain unstable homotopy sheaves of spheres. The following section, section \ref{section:connectivity_theorem}, discusses homotopy sheaves of spaces, characterizing important properties. It also states the unstable connectivity theorem. The last section describes some beautiful applications to the study of algebraic vector bundles.

Some things we do not do in this article include Voevodsky's groundbreaking work on the Bloch-Kato and Milnor
conjectures. A superb overview of this is given \cite{MorelVoevodskyproofMilnor1998}. We also do not deal with any
stable results, which have their own article in this handbook \cite{Isaksen-Ostvaer-handbook}.

In our presentation, we concentrate on those aspects and applications of the theory that relate to the calculation of the unstable homotopy
sheaves of spheres. The most notable of the applications, at present, is the formation of an $\Aone$ obstruction theory
of $B\GL_n$, and from there, the proof of strong results about the existence and classification of vector bundles on smooth
affine varieties. These can be found in Section \ref{sec:VB}.

\subsection{Acknowledgements} We wish to thank Aravind Asok, Jean Fasel, Raman Parimala and Joseph Rabinoff, as well as the organizers of Homotopy
Theory Summer Berlin 2018, the Newton Institute and the organizers of the Homotopy Harnessing Higher Structures
programme. The first-named author was partially supported by National Science Foundation Award DMS-1552730. 

\section{Overview of the construction of unstable $\Aone$-homotopy theory}\label{Section:construction_unstable_A1-homotopy_theory}

\subsection{Homotopy theory of Spaces} 

Morel and Voevodsky \cite{morelMathbbAHomotopyTheory1999} constructed a ``homotopy theory of schemes,'' in a category
sufficiently general to include both schemes and simplicial sets. By convention, the objects in such a category may be
called a ``motivic space'' or an ``$\Aone$-space'' or, most commonly, simply a ``space''. The underlying category of
spaces is a category of simplicial presheaves on a category of schemes, where schemes themselves are embedded by means
of a Yoneda functor. Two localizations are then performed: a Nisnevich localization and an $\Aone$-localization. Before
describing these, we clarify the notion of ``a homotopy theory'' that we use. Two standard choices are that a homotopy
theory is a simplicial model category or a homotopy theory is an $\infty$-category, a.k.a., a quasi-category. Background
on simplicial model categories can be found in \cite{HirschhornModelCategoriesTheir2003} and
\cite{Riehl_Categorical_homotopy_theory}. Lurie's {\em Higher Topos Theory} \cite{LurieHighertopostheory2009} contains
many tools used for doing $\Aone$-homotopy theory with $\infty$-categories. In the present case, either the
model-categorical or the $\infty$-categorical approach can be chosen, and once this choice has been made, presheaves of
simplicial sets and the resulting localizations take on two different, although compatible, meanings, either of which
produces a homotopy theory of schemes as described below.

To fix ideas, we use simplicial model categories. Let $S$ be a noetherian scheme of finite dimension, and let $\Sm_S$
denote the category of smooth schemes of finite type over $S$, defined as a full subcategory of schemes over $S$. Let
$\cat{sPre}(\Sm_S)$ denote the category of functors from the opposite category of $\Sm_S$ to simplicial sets, i.e.,
$\cat{sPre}(\Sm_S) = \Fun(\Sm_S^{\op},\sSet)$. Objects of $\cat{sPre}(\Sm_S)$ are called \textit{simplicial
 presheaves}. The properties of $\cat{sPre}(\Sm_S)$ that make it desirable for homotopy theoretic purposes are the following:
\begin{itemize}
\item There is a \textit{Yoneda embedding} $\eta: \Sm_S \to\cat{sPre}(\Sm_S)$, sending a scheme $X$ to the simplicial
 presheaf $U \mapsto \Sm_S(U,X)$, where the set of maps is understood as a $0$-dimensional simplicial set.
\item The category $\cat{sPre}(\Sm_S)$ has all small limits and colimits. In this it is different from the category $\Sm_S$ itself,
 which is far from containing all the colimits one might wish for in doing homotopy theory. The Yoneda embedding is
 continuous, but not cocontinuous.
\item The category $\cat{sPre}(\Sm_S)$ is simplicial.
\end{itemize}
The Yoneda embedding will be tacitly used throughout to identify a smooth scheme $X$ with the presheaf it represents, a $0$-dimensional
object of $\Spaces[S]$. We remark in passing that $S$ itself represents a terminal object of $\Spaces[S]$.

One can put several, Quillen equivalent, model structures on $\cat{sPre}(\Sm_S)$ such that weak equivalences are
detected objectwise, meaning that a map $X \to Y$ of simplicial presheaves is a weak equivalence if and only if the map
$X(U) \to Y(U)$ of simplicial sets is a weak equivalence for all $U$ in $\Sm_S$. Some standard choices for such a model
structure are the injective \cite{JardineSimplicialpresheaves1987}, projective \cite{blanderLocalProjectiveModel2001}
and flasque \cite{IsaksenFlasquemodelstructures2005} \cite{DuggerHypercoverssimplicialpresheaves2004}. We then carry out
two left Bousfield localizations. One can learn about Bousfield localization in
\cite{HirschhornModelCategoriesTheir2003}.

The first localization is analogous to the passage from presheaves to sheaves and depends on a choice of Grothendieck
topology on $\Sm_S$. We remind the reader that a Grothendieck topology is not a topology in the point--set sense, rather
it is sufficient data to allow one to speak meaningfully about locality and sheaves. We refer to
\cite{maclaneSheavesGeometryLogic1992} for generalities on sheaves.

The standard choice of Grothendieck topology for $\Aone$-homotopy theory is the {\em Nisnevich} topology, although the
\'etale topology is also used, producing a different homotopy theory of spaces. If one wishes to study possibly
nonsmooth schemes, then different topologies, for instance the cdh topology, can be used. The theory in the nonsmooth
case is less well developed, but \cite{voevodskyUnstableMotivicHomotopy2010} establishes that the theory for all schemes
and the cdh topology is a localization of the theory for smooth schemes and the Nisnevich topology, so that there is a
localization functor $L$ from one homotopy category to the other. Strikingly, the same
paper proves that 
that for a map $f :X \to Y$ in the smooth, Nisnevich theory, $Lf$ is a weak equivalence if and only if $\Sigma f$, the
suspension, is a weak equivalence.

In order to describe the Nisnevich topology, it is first necessary to describe \'etale maps: A map
$Y \to X$ in $\Sm_S$ is \textit{\'etale} if for every point $y$ of $Y$, the induced map on tangent spaces is an
isomorphism \cite[2.2 Proposition 8 and Corollary 10]{Bosch_Lutkebohmert_Raynaud-Neron_models}. Many other
characterizations of \'etale morphisms can be found in \cite[\S 17]{GrothendieckElementsgeometriealgebrique1967}
. A finite collection of maps
$\{ U_i \to X\}$ is a \textit{Nisnevich cover} if the maps $U_i \to X$ are \'etale and for every point $x \in X$, there
is a point $u$ in some $U_i$ mapping to $x$ such that the induced map on residue fields $\kappa(x) \to \kappa(u)$ is an
isomorphism. The topology generated by this pretopology is the Nisnevich topology on $\Sm_S$. We remark that the
additional constraint imposed by the Nisnevich condition is not vacuous even when $S = \Spec \bar k$ is the spectrum of
an algebraically closed field. For
instance, when $|n|>1$, the $n$-th power map $\C^\times \to \C^\times$---or more correctly
$\mathbb{G}_{m,\CC} \to \mathbb{G}_{m,\CC}$---is an \'etale cover but not a Nisnevich cover because the map on the residue fields of the
generic points is the homomorphism $\CC(x) \to \CC(x)$ sending $x \mapsto x^n$, which is not an isomorphism. One sees that
given Nisnevich cover $f: \coprod U_i \to X$, there must be a dense open subset $V_0$ of $X$ on which $f$ admits a
section, and a dense open subset $V_1$ of the complement of $V_0$ on which $f$ admits a (possibly different) section and
so on, so that $X$ has a stratification such that $f$ admits sections on the open complements of the strata; this
definition is used in \cite[Section 3.1]{deligneLecturesMotivicCohomology2001}.

A hypercover is a generalization of the \v{C}ech cover
\[\ldots \rightthreearrow \coprod_i (U_i \cap U_j) \rightrightarrows \coprod_i U_i \to X \] associated to the cover
$\{ U_i \to X\}$. See, for example \cite[Chapter 3]{FriedlanderEtalehomotopysimplicial1982a} or \cite[\S
4]{DuggerHypercoverssimplicialpresheaves2004}. One may therefore define hypercovers for the Nisnevich topology. The
first localization we carry out forces the maps $\hocolim_n U_n \to X$ to be weak equivalences for every hypercover
$U \to X$. This localization is called \textit{Nisnevich localization} when the corresponding topology is the Nisnevich
topology. The resulting simplicial model structures are referred to as the local model structures, and they are all
Quillen equivalent via identity functors. The local injective model structure on presheaves was originally constructed by Jardine in
\cite[Theorem 2.3]{JardineSimplicialpresheaves1987}, in a slight generalization of a construction due to Joyal
\cite{joyalLetterGrothendieck1983} who considered only sheaves rather than presheaves. What we call the ``local
injective'' model structure is notably referred to the ``simplicial model structure'' in
\cite{morelMathbbAHomotopyTheory1999}.

The theory constructed by Joyal and Jardine is very general, applying to all sites. We refer to
\cite{jardineLocalHomotopyTheory2015} for more about local homotopy theory \textit{per se}. In the cases called for
by $\Aone$-homotopy, there is a more elementary approach to the theory. Building on ideas of Brown \& Gersten
\cite{brownAlgebraicTheoryGeneralized1973}, one defines a \textit{distinguished square} of schemes to be a diagram
\[ \xymatrix{ U \times_X V \ar[r] \ar[d] & V \ar^p[d] \\ U \ar^i[r] & X } \] of schemes where $i$ is an open immersion,
$p$ is an \'etale morphism and $p$ restricts to an isomorphism on the closed complements $p^{-1} ( X - U) \to X - U$
(given the reduced induced subscheme structures). Then a functor $\sh F: \Sm_S^\op \to \Set$ is a Nisnevich sheaf if
$\sh F$ takes distinguished squares to cartesian squares of sets, and a functor $\sh X: \Sm_S^\op \to \sSet$ is, loosely
speaking, suitable for $\Aone$-homotopy theory if it takes distinguished squares to homotopy cartesian diagrams of
simplicial sets. The following result, \cite[Lemma 4.1]{blanderLocalProjectiveModel2001}is the simplest in a family of
many such results.
\begin{pr}
 A simplicial presheaf $\sh X : \Sm_S^\op \to \sSet$ is fibrant in the projective local model structure if $\sh X(\cdot)$
 takes values in Kan complexes and $\sh X$ takes distinguished squares to homotopy cartesian squares.
\end{pr}
Variations on this idea are considered in \cite[Section 4.2]{IsaksenFlasquemodelstructures2005}, \cite[Section
3.2]{asokAffineRepresentabilityResults2017} and \cite{Voevodsky-Homotopy_theory_simplicial_sheaves_cd_topologies}, where
we have listed the references in order from least to most general.

The second localization makes the projection maps $X \times \Aone^1 \to X$ into weak equivalences for all smooth
schemes $X$ in $\Sm_S$. Any resulting simplicial model structure is an \textit{$\Aone$-model structure}. A pleasant
overview of the relationships between the various model categories appearing in this story can be found in
\cite{IsaksenFlasquemodelstructures2005}.

\begin{definition}
Let $R_{\Aone}$ denote a fibrant replacement functor in the injective $\Aone$ model structure. We will say that an object $\sh
X$ is \textit{$\Aone$-local} if the map $\sh X \to R_{\Aone} \sh X$ is a weak equivalence in a local model structure. 
\end{definition}
That is
$\Aone$-local objects are the objects for which the local homotopy type already recovers the $\Aone$-homotopy type.

These simplicial model categories of spaces allows us to carry out homotopy theory on schemes. For example, we can form limits
and colimits, in particular smash products of schemes or of schemes and simplicial sets, and give meaning to the spaces
in \eqref{deghom}. Many results from the classical homotopy theory of simplicial sets carry over. For example, we have
excision: suppose that $Z \to X$ is a closed immersion in $\Sm_k$ 
and $Y \to X$ is an open subset of $X$. Then there is a pushout square of
schemes \begin{equation}\label{CoverofXwithYandZ} \xymatrix{ Y - Z \ar[r] \ar[d] & X-Z \ar[d] \\ Y \ar[r]&
 X} \end{equation} corresponding to a Nisnevich cover of $X$. The Nisnevich localization procedure guarantees that
Nisnevich covers give rise to pushouts in spaces, causing \eqref{CoverofXwithYandZ} to be a homotopy pushout, giving the
excision weak-equivalence \[ X/(X-Z) \weq Y/(Y-Z)
.\]

We also can take any vector bundle $p:E \to X$, where $X$ is a smooth scheme, and decompose $X$ as an open cover of
subschemes $U_i$ such that the induced maps $p^{-1}(U_i) \to U_i$ is isomorphic to the projection
$\A^n \times U_i \to U_i$, which is a weak equivalence. A colimiting argument then shows that $p: E \to X$ is itself an
$\Aone$-weak equivalence. This works for any sort of Nisnevich-locally-trivial fibration with
$\Aone$-contractible fibres, and in particular, for quasiprojective varieties over a field $k$, the Jouanolou trick \cite[Proposition
4.3]{WeibelHomotopyalgebraictheory1989} produces a map $p: \Spec R \to X$ which is an $\Aone$-equivalence having affine
source. Any quasiprojective smooth $k$-variety is therefore $\Aone$-equivalent to an affine variety.

\subsection{Homotopy Sheaves}

An inconvenience in $\Aone$-homotopy theory is varying availability of basepoints. By a ``basepoint'' we might mean a
map $p: S \to \sh X$ over $S$, viz., a morphism from the terminal object of the category $\Spaces[S]$. Unfortunately, since $\sh X$
amounts to a family of simplicial sets, parametrized by $\Sm_S^\op$, for certain objects $U \to S$ of $\Sm_S^\op$, there
may be path components of $\sh X(U)$ that are not in the image of the map $\pi_0\sh X(S) \to \pi_0 \sh X(U)$. For
instance, if $V$ is the closed subscheme of $\A^2$ over $\RR$ determined by the vanishing of $x^2 + y^2 + 1$, then
$V(\RR)$ is empty, while $V(\CC)$ is a discrete set of uncountably many points.

In order to handle this technicality, one must allow ``basepoints'' $x_0 \in \sh X(U)$ that are not in the image
of $\sh X(S) \to \sh X(U)$. That is, for any object $U$ of $\Sm_S$, and any $x_0 \in (\sh X(U))_0$ one defines a
sequence of presheaves on the slice category $\Sm_S/U$ by
\begin{equation} \label{eq:pipre} (V \overset{f}{\to} U) \mapsto \pi_n( |\sh X(V)|, f^*(x_0)). \end{equation}
Here $|\sh X(V)|$ denotes the geometric realization of the simplicial set.

The Nisnevich sheaves associated to the functors of \eqref{eq:pipre} are the \textit{homotopy sheaves} of $\sh X$ for
the basepoint $x_0$, and will be denoted $\bpi_n(\sh X, x_0)$. If $n \ge 1$, then the sheaves are sheaves of groups, and
they are abelian if $n\ge 2$. One may define an unpointed $\bpi_0(\sh X)$ similarly.

We define $\bpi_n^{\Aone}$, the \textit{$\Aone$-homotopy sheaves}, by first replacing $\sh X$ by an $\Aone$ local object---for
instance, $R_{\Aone} \sh X$, then calculating $\bpi_n$ of the resulting object. Both the homotopy sheaves and the $\Aone$-homotopy sheaves satisfy a
Whitehead theorem (in fact, this is taken as the definition of local weak equivalence in \cite{JardineSimplicialpresheaves1987}).
\begin{pr}
 If $f: \sh X \to \sh Y$ is a morphism in $\Spaces[S]$, then $f$ is an $\Aone$-weak equivalence if and only if the morphisms
 \begin{equation*}
 f_* : \bpi^{\Aone}_0(\sh X) \to \bpi^\Aone_0(\sh F)
 \end{equation*}
 and 
 \begin{equation*}
 f_* :\bpi^\Aone_n(\sh X, x_0) \to \bpi^\Aone_n(\sh Y, f(x_0))
 \end{equation*}
 are isomorphisms for all choices of basepoint $x_0 \in (\sh X(U))_0$. 
\end{pr}

The story is much simpler in the case where $\sh X$ is $\Aone$-connected, i.e., when $\bpi_0^\Aone (\sh X)$ is a singleton. In
this case, non-globally-defined basepoints need not be considered. This version appears as \cite[Proposition 2.14 , \S
3]{morelMathbbAHomotopyTheory1999}.

One can define homotopy groups of $\sh X$ by taking sections of the homotopy sheaves. One has the following
\begin{pr}
 Let $k$ be a field and let $\sh X$ be an object of $\Spaces$. Choose a basepoint $x_0 \in (\sh X (k))_0$. Then
 \[ \bpi_n^\Aone(\sh X, x_0)(k) = [ S^n \wedge (\Spec k)_+, \sh X]_\Aone \]
\end{pr}

\subsection{Spheres}

Let $\Gm$ denote the punctured affine line $\bbb{G}_m=\Aone - \{0 \}$, given the basepoint $1$, and let $S^1$ denote the
pointed simplicial circle. The role of the spheres $S^n$ in classical algebraic topology is now played by smash products
of the spaces $S^1$ and $\Gm$. We use the notation
$$S^{p + q \alpha} = S^{p+q,q} := (S^1)^{\wedge p} \wedge \G_m^{\wedge q},$$ and such spaces will be called
\textit{spheres}. The presence of two different indexing conventions for spheres in $\Aone$-homotopy theory can be
confusing, but both appear in the literature, so we give both as well. The notation $\Sigma X$ (or $\Sigma_{S^1} X$ if
there is possible confusion) denotes $X \wedge S^1$. The question of which schemes have the $\Aone$-homotopy type of
spheres has been studied in \cite{AsokSmoothmodelsmotivic2017}, and the most common examples are the following.

\begin{example}
The scheme $\bbb{P}^1$ is $\Aone$-homotopy equivalent to $S^{1+\alpha} = S^{2,1}$. Specifically, the pushout diagram \[\xymatrix{\G_m \ar[r] \ar[d] & \Aone \ar[d] \\ 
\Aone \ar[r] & \bbb{P}^1}\] and the weak equivalence $\Aone \weq \ast$, induce a weak equivalence $$\bbb{P}^1 \simeq \Sigma \G_m.$$
\end{example}

\begin{example}\label{An-0isS2n-1n}
The scheme $\bbb{A}^n - \{0\}$ is $\Aone$-homotopy equivalent to $S^{n-1+n\alpha} = S^{2n-1,n}$. The case $n=1$ is immediate, and then one may proceed by induction on $n$: the pushout diagram \[\xymatrix{\G_m \times (\bbb{A}^{n-1} - \{0\}) \ar[r] \ar[d] & \Aone\times (\bbb{A}^{n-1} - \{0\}) \ar[d] \\ 
\G_m \times \bbb{A}^n \ar[r] & \bbb{A}^n - \{0\}}\] and the weak equivalences $\G_m \times \bbb{A}^n \weq \G_m$ and $\Aone\times (\bbb{A}^{n-1} - \{0\}) \weq \bbb{A}^{n-1} - \{0\}$, induce a weak equivalence between $\bbb{A}^n - \{0\}$ and the homotopy pushout of the diagram $Y \leftarrow X \times Y \rightarrow X$ of projection maps for $X = \G_m$ and $Y=\bbb{A}^{n-1} - \{0\}$. Since this latter homotopy pushout is identified with $\Sigma (X \wedge Y)$ by pushing out the rows of the diagram $$ \xymatrix{ X \ar[d] & \ar[l] X \vee Y \ar[d] \ar[r] & Y \ar[d]\\ 
X \ar[d] &\ar[l] X \times Y \ar[d] \ar[r] & Y \ar[d]\\ \ast & \ar[l] X \wedge Y \ar[r] & \ast,}$$ the result follows. 
\end{example}

An important related example is the identification of the $\Aone$-homotopy type of $\bbb{P}^n_k/\bbb{P}^{n-1}_k$ with
$S^{n+n\alpha} = S^{2n,n}$.

\begin{example}\label{Pn/Pn-1=Sn+nalpha}
There is a standard closed immersion\footnote{We use the algebre-geometric term ``closed immersion'' for a map
 isomorphic to the inclusion of a closed subscheme. The usual term in topology is ``closed embedding'', which is also
 used in \cite{morelMathbbAHomotopyTheory1999}, but is not widespread in the literature on $\Aone$-homotopy theory.} $\bbb{P}^{n-1} \hookrightarrow \bbb{P}^n$ sending homogeneous coordinates
$[x_0: \ldots: x_n]$ to $[0;x_0: \ldots:x_n]$. The homotopy type of the quotient space $\bbb{P}^n/\bbb{P}^{n-1}$ is
identified by the chain of weak equivalences $$\bbb{P}^n/\bbb{P}^{n-1} \weq \bbb{P}^n/(\bbb{P}^n -\{0\}) \weq \bbb{A}^n/(\bbb{A}^n -\{0\}) \weq \Sigma(\bbb{A}^n -\{0\})$$ with $\Sigma(\bbb{A}^n -\{0\})$. Therefore, by example \ref{An-0isS2n-1n}, we have $\bbb{P}^n/\bbb{P}^{n-1} \weq S^{n+n\alpha}$.
\end{example}

\subsection{Thom spaces and purity}

Let $V \to X$ be a vector bundle and let $X \hookrightarrow V$ be the zero section. 

\begin{df}
The {\em Thom space}, denoted $\Th(V) $ or $X^V$, of $V$ is defined $$\Th(V) := V/(V -X).$$ 
\end{df}

A parametrized version of Example \ref{Pn/Pn-1=Sn+nalpha} shows that the $\Aone$-homotopy type of the Thom space over a scheme can be alternatively described as the cofiber of a closed immersion. Namely, let $1$ denote the trivial bundle on $X$. For any vector bundle $V$ on $X$, let $\bbb{P}(V)$ denote the projective space bundle given by the fiberwise projectivization. As in Example \ref{Pn/Pn-1=Sn+nalpha}, there is a standard closed immersion $\bbb{P}(V) \hookrightarrow \bbb{P}(V \oplus 1)$. The complementary open subscheme is isomorphic to $V$.

\begin{pr}
There is an $\Aone$-weak equivalence $\Th(V) \simeq \bbb{P}(V \oplus 1)/\bbb{P}(V)$.
\end{pr}

\begin{proof}
By excision, there is an $\Aone$-weak equivalence $\Th(V) \simeq \bbb{P}(V \oplus 1)/ (\bbb{P}(V \oplus 1) - X)$. The scheme $\bbb{P}(V \oplus 1) - X$ is the total space of the vector bundle $\mathcal{O}(1)$ over $\bbb{P}(V)$, giving the claimed weak equivalence.
\end{proof}

The analogue of the tubular neighborhood theorem from classical topology is Morel--Voevodsky's Purity
theorem.

\begin{namedthm}{Purity Theorem}[Morel, Voevodsky]\label{Puritythm}
 Let $Z \hookrightarrow X$ be a closed immersion in $\Sm_S$. Then there is a natural $\Aone$-weak
 equivalence $$ \Th(N_Z X) \simeq X/(X-Z) $$ where $N_Z X \to Z$ denotes the normal bundle.
\end{namedthm}

The original proof is \cite[Section 3 Theorem 2.23 p.115]{morelMathbbAHomotopyTheory1999}, and a particularly accessible
exposition, done over an algebraically closed field, may be found in \cite[Section 7]{AntieauElmanto-primer}, using
unpublished notes of A.~Asok and \cite{Hoyoissixoperationsequivariant2017}.

Here, we will give an outline of a proof and draw some pictures in the case
where $Z \hookrightarrow X$ is the inclusion of the point into the affine line $\{0\} \hookrightarrow \A^1_k$.

First, we explain the concept of ``blow up'' in algebraic geometry briefly. For a detailed treatment, we refer to
\cite{hartshorneAlgebraicGeometry1977}. Given a closed immersion $Z \hookrightarrow X$, the blow-up $\pi: \Bl_Z X \to X$
of $X$ along $Z$ is a map satisfying the property that the restriction of $\pi$ to $\Bl_Z X - \pi^{-1}(X)$ is an
isomorphism and the inverse image of $Z$, called {\em the exceptional divisor}, is $\pi^{-1} Z \cong \bbb{P} N_Z X$ the
projectivization of the normal bundle. In other words, we cut out $Z$ and glue in a projectivization of the normal
bundle. In the topology of manifolds, this can be accomplished by first removing an open tubular neighborhood of $Z$
from $X$, so that one has introduced a boundary component $S$ homeomorphic to a sphere bundle in the normal bundle
$N_Z X$, and then taking a fibrewise quotient of $S$ to replace it by a projectivization of $N_z X$. Algebraically, one
takes the sheaf of ideals $\mathcal{I}$ determining the closed immersion $Z \to X$, forms the sheaf of graded algebras
$\oplus_{i=0}^{\infty}\mathcal{I}^i$ and takes the relative $\Sheafproj$ construction, i.e., $\Bl_Z X \to X$ is the
canonical map $\Sheafproj_X \oplus_{i=0}^{\infty}\mathcal{I}^i \to X$. The morphism $\pi:\Bl_Z X \to X$ satisfies the
property that the sheaf of ideals $(\pi^{-1}\mathcal{I} )\mathcal{O}_{\Bl_Z X}$ associated to $\pi^{-1}Z$ is the sheaf
of sections of a line bundle, and any other morphism $Y \to X$ with this property factors through $\pi$, \cite[II
7.14]{hartshorneAlgebraicGeometry1977}.

\begin{figure}[ht!]
\centering
\includegraphics[angle=0,origin=c,width=90mm]{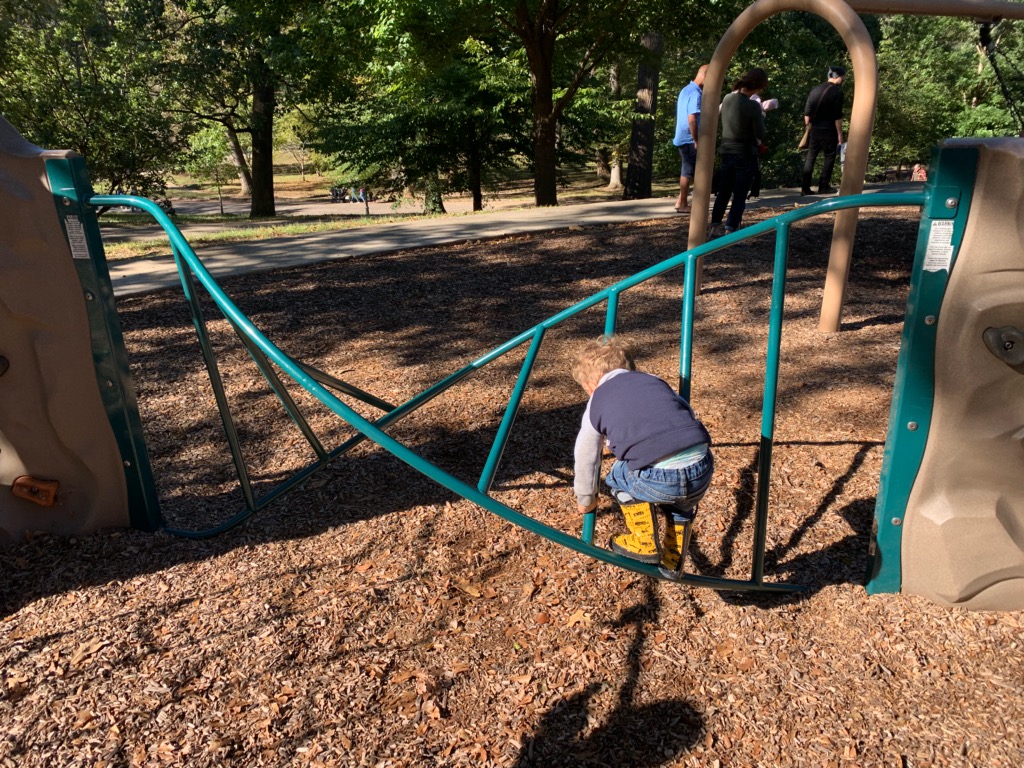}
\caption{$\Bl_0 \bbb{A}^2$ Jungle Gym \label{blowup picture}}
\end{figure}

\begin{example} 
 The blow-up $\Bl_0 \bbb{A}^n$ of affine $n$-space $\bbb{A}^n$ at the origin is the closed subscheme of
 $\bbb{A}^n \times \bbb{P}^{n-1}$ determined by $\langle x_i y_j -x_j y_i \rangle$ where $(x_1, \ldots, x_n)$ are the
 coordinates on $\bbb{A}^n$ and $[y_1, \ldots, y_n]$ are homogeneous coordinates on $\bbb{P}^n$. The restriction to
 $\Bl_0 \bbb{A}^n$ of the projection $\bbb{A}^n \times \bbb{P}^{n-1} \to \bbb{P}^{n-1}$ is the projection of the
 tautological bundle, denoted $\calO(-1)$, on $\bbb{P}^{n-1}$ to $\bbb{P}^{n-1}$. See Figure
 \ref{Figure_for_Blowup_example} for the case $n=2$. The green line is the exceptional divisor $\bbb{P}^1$, and the red
 and blue lines are both fibers of the line bundle $\calO(-1)$. Blow-ups of $\bbb{A}^2$ are sometimes used as jungle
 gyms. See Figure \ref{blowup picture}. \end{example}

\begin{figure}
\begin{tikzpicture}[myxyz,scale=1.5]
 \foreach \a in {36,42,...,215} {
 \fill[black!10!white, opacity=.7]
 let \n1 = {2*cos(\a)}, \n2 = {2*sin(\a)}, \n5 = {\a+6},
 \n3 = {2*cos(\n5)}, \n4 = {2*sin(\n5)} in
 ( \n1, \n2, 0.02*\a ) -- (-\n1, -\n2, 0.02*\a ) --
 (-\n3, -\n4, 0.02*\n5) -- ( \n3, \n4, 0.02*\n5) -- cycle;
 \draw[help lines]
 let \n1 = {2*cos(\a)}, \n2 = {2*sin(\a)}, \n5 = {\a+6},
 \n3 = {2*cos(\n5)}, \n4 = {2*sin(\n5)} in
 ( \n1, \n2, 0.02*\a ) -- (-\n1, -\n2, 0.02*\a )
 ( \n1, \n2, 0.02*\a ) -- ( \n3, \n4, 0.02*\n5)
 (-\n1, -\n2, 0.02*\a ) -- (-\n3, -\n4, 0.02*\n5);
 }

 \draw[very thick, green!50!black] (0,0,.02*36) -- (0,0,.02*216);
 \draw[very thick, blue!70!black] (-2, 0, .02*180) -- (2, 0, .02*180);
 \draw[very thick, red!80!black] (0, -2, .02*90 ) -- (0, 2, .02*90);

 \node[left] at (-2.3, 0, .02*140) {$\Bl_0(\A^2)$};
 \node[right, red!80!black, font=\small] at (0.1, 0, .02*90) {line slope $\infty$};
 \node[right, blue!70!black, font=\small] at (2.1, 0, .02*180) {line slope $0$};
 \node[above right, green!50!black, font=\small, inner sep=0]
 at (0, 0, .02*216) {$\pi^{-1}(0)$};

 \draw[very thick, ->] (0, 0, -.2) to["$\pi$"] (0, 0, -1.2);

 \begin{scope}[yshift=-3cm]
 \begin{scope}[transformxy]
 \fill[black!10!white, opacity=.7] (-2, -2) rectangle (2, 2);
 \draw[help lines] (0,0) circle[radius=2];
 \draw[help lines] foreach \a in {36,42,...,215} {
 ($-2*({cos(\a)}, {sin(\a)})$) -- ($2*({cos(\a)}, {sin(\a)})$)
 };
 \draw[very thick, blue!70!black] (-2, 0) -- (2, 0);
 \draw[very thick, red!80!black] (0, -2 ) -- (0, 2);
 \point[green!50!black] at (0,0);
 
 \end{scope}
 \node[left] at (-2.3, 0) {$\A^2$};
 \end{scope}

\end{tikzpicture}
\begin{caption}\relax
Blow-up 
\label{Figure_for_Blowup_example}
\end{caption}
\end{figure}

An essential component of the proof of the Purity Theorem \ref{Puritythm} is the ``deformation to the normal bundle'', a
pre-existing idea in intersection theory\cite[Chapter 5]{FultonIntersectionTheory1984}. The input to this construction
is a closed immersion $Z \hookrightarrow X$ in $\Sm_k$ and the output is a family of closed immersions over
$\Aone$: \begin{equation}\label{DeformationToNormalCone}\xymatrix{ Z \times \Aone \ar[rd] \ar[rr] && D_Z X \ar[ld] \\ &
 \Aone & } \end{equation} such that over $t=1$, the fiber of the family is the original closed immersion
$Z \hookrightarrow X$, and over $t=0$, the fiber of the family inclusion of the zero section into $N_Z X$.

To form this family, first consider $\Bl_{Z \times 0} (X \times \Aone) \stackrel{\pi}{\to} X \times \Aone \to
\Aone$. The inverse image of $t=0$ under the composition is $\pi^{-1}(X \times \{0\})$. It is the union of the
exceptional divisor $\bbb{P}(N_{Z \times 0} (X \times \Aone))$ and a copy of $\Bl_{Z \times 0} (X \times 0)$. These are
glued by identifying the copy of $\bbb{P} N_{Z}X$ embedded in $\Bl_{Z \times 0} (X \times 0)$ as the exceptional
divisor, and the following copy of $\bbb{P} N_{Z}X$ in $\bbb{P}(N_{Z \times 0} (X \times \Aone))$: the bundle
$N_{Z \times 0} (X \times \Aone)$ is canonically isomorphic to $N_Z X \oplus 1$. Thus
$\bbb{P}(N_{Z \times 0} (X \times \Aone))$ is the union of $N_{Z} X$ and $\bbb{P} N_Z X$, and the latter $\bbb{P} N_Z X$
is the copy we seek. Then let $$D_Z X : = \Bl_{Z \times 0} (X \times \Aone) - \Bl_{Z \times 0} (X \times 0).$$ Since
the copy of $Z \times \Aone$-in $\Bl_{Z \times 0} (X \times \Aone)$ intersects the exceptional divisor as the zero
section of $N_{Z} X $ in $\bbb{P}(N_{Z \times 0} (X \times \Aone))$, the scheme $D_Z X$ provides the claimed
family. $D_0 \Aone$ is shown in Figure \ref{DeformationToNormalConePicture}.

\begin{figure}
\begin{tikzpicture}[myxyz]
 \foreach \a in {36,42,...,215} {
 \fill[black!10!white, opacity=.7]
 let \n1 = {2*cos(\a)}, \n2 = {2*sin(\a)}, \n5 = {\a+6},
 \n3 = {2*cos(\n5)}, \n4 = {2*sin(\n5)} in
 ( \n1, \n2, 0.02*\a ) -- (-\n1, -\n2, 0.02*\a ) --
 (-\n3, -\n4, 0.02*\n5) -- ( \n3, \n4, 0.02*\n5) -- cycle;
 \draw[help lines]
 let \n1 = {2*cos(\a)}, \n2 = {2*sin(\a)}, \n5 = {\a+6},
 \n3 = {2*cos(\n5)}, \n4 = {2*sin(\n5)} in
 ( \n1, \n2, 0.02*\a ) -- (-\n1, -\n2, 0.02*\a )
 ( \n1, \n2, 0.02*\a ) -- ( \n3, \n4, 0.02*\n5)
 (-\n1, -\n2, 0.02*\a ) -- (-\n3, -\n4, 0.02*\n5);
 }

 \draw[very thick, green!50!black] (0,0,.02*36) -- (0,0,.02*216);
 \draw[very thick, blue!70!black] (-2, 0, .02*180) -- (2, 0, .02*180);
 \draw[very thick, red!80!black] (0, -2, .02*90 ) -- (0, 2, .02*90);
 
 \draw[very thick, ->] (0, 0, -.2) -- (0, 0, -1.2);

 \node (L1) at (-4.8, 0, .02*140) {$\Bl_{Z\times 0}(X\times\A^1)$};

 \node[right, red!80!black, font=\small] at (0.1, 0, .02*90)
 {$\Bl_{Z\times0}(X\times 0)$};
 \node[above right, green!50!black, font=\small, inner sep=-2mm, align=center]
 at (0, 0, .02*216) {exceptional divisor\\$\mathbb{P}(N_z\oplus 1)$};

 \begin{scope}[yshift=-3.6cm]
 \begin{scope}[transformxy]
 \fill[black!10!white, opacity=.7] (-2, -2) rectangle (2, 2);
 \draw[help lines] (-2, -2) rectangle (2, 2);
 \draw[help lines] foreach \a in {36,42,...,215} {
 ($-2*({cos(\a)}, {sin(\a)})$) -- ($2*({cos(\a)}, {sin(\a)})$)
 };
 \draw[very thick, blue!70!black] (-2, 0) -- (2, 0);
 \draw[very thick, red!80!black] (0, -2 ) -- (0, 2);
 \point[green!50!black] at (0,0);

x \node (A) at (0, -2.4) {$\A^1$};
 \draw[->] (A) -- (-2,-2.4);
 \draw[->] (A) -- (2,-2.4);

 \node (X) at (2.4, 0) {$X$};
 \draw[->] (X) -- (2.4,-2);
 \draw[->] (X) -- (2.4,2);

 \end{scope}

 \node (L2) at (-4.8, 0) {\llap{\color{blue!70!black}$Z\times\A^1\to$} $X\times\A^1$};

 \draw[very thick, ->] (0, 0, -2) -- (0, 0, -3);
 \end{scope}

 \begin{scope}[yshift=-7.1cm]
 \begin{scope}[transformxy]
 \draw[help lines, thick] (-2, 0) -- (2, 0);
 \point["$t=0$" below, red!80!black] at (0, 0);

 \node (L3) at (-4.8, 0) {$\A^1$};

 \end{scope}
 \end{scope}

 \draw[thick, ->] (L1) -- (L2);
 \draw[thick, ->] (L2) -- (L3);
 \draw[thick, ->] (L1) to[bend left, "$f$"] (L3);

 \node at (-2.4, 0, -8.5) {$D_z(X) = \Bl_{Z\times 0}(X\times\A^1) - \color{red!80!black}\Bl_{Z\times 0}(X\times 0)$};

\end{tikzpicture}
\begin{caption}\relax
Deformation to the normal bundle
\label{DeformationToNormalConePicture}
\end{caption}
\end{figure}

To prove purity, one shows that the fibers of \eqref{DeformationToNormalCone} above $t=0$ and $t=1$ induce $\Aone$-weak
equivalences $$i_0 : N_Z X / (N_Z X - Z)\to D_Z X/ (D_Z X - Z \times \Aone_k)$$ and $$i_1 : X/(X-Z) \to D_Z X/ (D_Z X - Z
\times \Aone),$$ respectively. Like the \'etale topology, the Nisnevich topology is fine enough to reduce many
arguments about closed immersions in $\Sm_S$ to arguments for the standard inclusion $\bbb{A}^n \hookrightarrow
\bbb{A}^{n+c}$. A precise theorem enabling such reductions may be found in \cite[II 4.10]{GrothendieckRevetementsetalesgroupe1971}. 

In the case where $Z \hookrightarrow X$ is $0 \hookrightarrow \bbb{A}^n$, the deformation to the normal bundle
$$D_Z X = \Bl_0 (\bbb{A}^n \times \Aone) - \Bl_0 (\bbb{A}^n \times 0) \cong \calO_{\bbb{P}^n_k} (-1) -
\calO_{\bbb{P}^{n-1}_k} (-1) \cong \calO_{\bbb{P}^n_k} (-1) \vert_{\bbb{A}^n_k} $$
is the restriction of the total space of the tautological bundle on $\bbb{P}^n_k$ restricted to the standard copy of
$\bbb{A}^n$. The fiber above $0$ is the zero-section of this bundle. The fiber above $1$ also defines a section. We can
see directly that $i_0$ and $i_1$ are $\Aone$-weak equivalences.

\section{Realizations}\label{RealizationsSection}

\subsection{Complex Realization}

 Suppose $k \subset \CC$ is a subfield of the complex numbers, then one may view the
 category $\Sm_k$ of smooth $k$-schemes as a category of complex varieties. Let us abuse notation and write $X(\CC)$
 for the complex manifold determined by a smooth $k$-scheme $\Sm_k$. The object $X(\CC)$, being a manifold, is the sort
 of thing to which classical homotopy theory is well adapted. One therefore might hope for a functor
 $|\cdot|: \Spaces \to \cat{Top}$ with good homotopical properties. Such a functor can be constructed in the most obvious way:
 for a smooth $k$-scheme $X$, one defines $|X| = X(\CC)$, as discussed above. The functor
 \[ |\cdot| : \Spaces \to \Top \]
 is then extended to all simplicial presheaves in such a way as to preserve (homotopy) colimits. For formal reasons,
 if $K$ is a simplicial set, viewed as a constant simplicial presheaf in $\Spaces$, then this realization $|K|$ agrees
 with the usual topological realization of the simplicial set $K$.

 The technical details in showing that this realization procedure works out is to be found in
 \cite{DuggerUniversalHomotopyTheories2001}, \cite{DuggerHypercoverssimplicialpresheaves2004} and
 \cite{DuggerTopologicalhypercoversmathbbA2004}. The first paper shows that the motivic projective model structure on $\Spaces$ is an $\Aone$-localization
 of a ``universal'' model structure, and then the last two show that the relations imposed in the motivic projective model
 structure, namely making hypercovers and trivial line bundles into weak equivalences, are relations that yield
 ordinary homotopy equivalences after realization. The upshot is a Quillen adjunction
 \[ | \cdot | : \Spaces \leftrightarrows \Top: S, \]
 when the left hand side is endowed with the motivic projective model structure. The approach outlined above is an
 improvement on a result of \cite{morelMathbbAHomotopyTheory1999} where a functor $t^\CC : \sh H(\CC) \to \sh H$ is
 produced, but only on the level of homotopy categories.

 Being a left adjoint, the realization functor constructed in this way is not compatible with fiber sequences, but one
 sees that it is compatible with products, in that $|X \times Y| \homeo |X| \times |Y|$,
 \cite[Section A.4]{PaninVoevodskyAlgebraicKTheory2009} and as a consequence, one may deduce
 \[ | B(X, G, Y)| \homeo B(|X|, |G|, |Y|) \]
 for two-sided bar constructions. This may suffice in certain cases.

 Since realization is homotopically well behaved, one obtains maps $[X, Y]_{\A^1_\CC} \to [|X|, |Y|]$. This map is
 understood when $X= S^{n+m \alpha}$ and $Y = S^{n+j\alpha}$, assuming $n \ge 2$ and $j \ge 1$, in which case one has maps
 \begin{equation} \label{eq:CComparison} \KMW_{j-m}(\CC) \to \pi_{n+m}(S^{n+j}). \end{equation}
 For dimensional reasons, this map vanishes when $j > m$. When $j=m$, the identity map generates both source and
 target, and one has an isomorphism $\ZZ = \GW(\CC) \to \ZZ = \pi_{n+m}(S^{n+m})$. In the cases where $j < m$, one
 knows from \cite{MorelmathbbAalgebraictopology2012} that $\KMW_{j-m}(\CC) = W(\CC)$ contains a single nontrivial
 class, that of $\eta^{j-m}$, the iterated Hopf map, and so the image of the realization map \eqref{eq:CComparison} is
 the group generated by $\eta_{\mathrm{top}}^{j-m}$, which vanishes if $j-m$ is sufficiently large.

\subsection{Real Realization}

Suppose now that $k = \RR$. Take $X(\CC) = \Mor_{\RR}(\Spec \CC , X)$ and recall that the Galois action
$C_2 = \Gal(\CC/\RR)$ endows this manifold with a $C_2$ action. Extending the functor $X \mapsto X(\CC)$ in such a way
as to preserve homotopy colimits, one produces an adjunction:
\[ | \cdot |_{\equi} : \Spaces \leftrightarrows C_2-\Top : S_{\RR}. \]
A map $f: X \to Y$ of $C_2$-spaces is a weak equivalence if and only if it induces weak equivalences after taking $H$
fixed points $f^H : X^H \to Y^H$ for both $H = C_2$ and $H=\{e\}$. With this definition, the functor $|\cdot |_{\equi}$ is
compatible with the homotopy---specifically \cite[Section 5.3]{DuggerTopologicalhypercoversmathbbA2004} establishes that
$|\cdot|_{\equi}$ is a left Quillen functor when the source is endowed with the motivic projective model structure, the
arguments being largely the same as in the complex case.

Composing $|\cdot |_{\equi}$ with the functor taking $C_2$ fixed points, which is a right adjoint, one arrives at a
functor $X \mapsto |X|_{\equi}^{C_2}:=|X|_{\RR}$, and if $X$ is a scheme, then $|X|_\RR$ is the analytic space $X(\RR)$
endowed with the obvious topology. The functor $|\cdot |_{\RR}$ cannot be a left adjoint, since it involves the taking
of fixed points, it does not preserve all colimits, but it does preserve smash products. Notably,
$| S^{n+m\alpha}|_\RR \weq S^n \wedge |\Gm^{\wedge m}|_\RR \weq S^n$, and so the induced maps
\[ \pi_{n+m\alpha} (S^{n+j\alpha})(\RR) = \KMW_{j-n}(\RR) \to \pi_n(S^n) = \ZZ \]
are quite different in character from those in the complex case of \eqref{eq:CComparison} above. The ring
$\KMW_\ast(\RR)$ is generated by classes $[a] \in \KMW_1$ where $a \in \RR^\times$ and $\eta \in \KMW_{-1}$. Under
$\R$-realization, $[a] \mapsto 0$ if $a > 0$, and $[a] \mapsto 1$ if $a < 0$, and $\eta \mapsto -2$, for a proof, see
\cite{AsokMotivicspheresimage2018a}[Proposition 3.1.3].

\subsection{\'Etale realization}

In the case of an arbitrary Noetherian base scheme, $S$, there is an \textit{\'etale realization} functor, constructed
in \cite{IsaksenEtalerealizationmathbbA2004}, operating along similar lines to the above two functors. In this case,
one begins with a functor $\Sm_S \to \text{pro-}\Spaces$, taking a smooth $S$-scheme $X$ to the \'etale topological type
of $X$, as defined in \cite{FriedlanderEtalehomotopysimplicial1982a}. One then extends this functor to all of
$\Spaces[S]$ by requiring it to commute with homotopy colimits. The result is again a left Quillen functor
\begin{equation*}
 | \cdot |_{\et} : \Spaces[S] \leftrightarrows \text{pro-}\Top : S 
\end{equation*}
with a certain extra complication since the target category is a category of pro-objects.

It might be hoped that if $S= \Spec k$ is a field, then $\Gal(\bar k/k)$ should act on $|X|_{\et}$ and the target
category might be enriched to $\Gal(\bar k/ k)$-equivariant-pro-spaces. The extent to which this can be done is not
known, owing in part to the technical difficulty of working with equivariant pro-spaces, and in part because it is known
that not all desirable properties for such a functor are attainable: \cite{Kassetalerealizationwhich2018}.

\section{Degree}\label{degree_section}

\subsection{The Grothendieck--Witt group}

Equipped with definitions of spheres and $\bbb{A}^1$-homotopy classes of maps, we consider again Morel's
$\bbb{A}^1$-degree homomorphism $$\deg: [S^n \wedge \Gm^{\wedge n}, S^n \wedge \Gm^{\wedge n}] \to \GW(k)$$ where the target is the Grothendieck--Witt group.

The Grothendieck--Witt group $\GW(k)$ is both complicated enough to support interesting invariants and simple
enough to allow explicit computations. It is defined to be the group completion of the semi-ring of isomorphism classes of $k$-valued, non-degenerate, symmetric, bilinear forms on finite dimensional $k$-vector spaces. In more detail: The isomorphism
classes of non-degenerate, symmetric, bilinear forms $\beta : V \times V \to k$ over $k$, where $V$ is a finite
dimensional $k$-vector space, admit the operations of perpendicular direct sum $\oplus$ and tensor product $\otimes$.
These operations give the set of such isomorphism classes the structure of a semi-ring. Taking the group completion, i.e., introducing formal differences of isomorphism classes, defines the ring $\GW(k)$.

$\GW(k)$ has a presentation \cite[II Theorem 4.1 pg 39]{LamIntroductionquadraticforms2005} given by generators $\langle a \rangle$ for $a$ in $k^\times$ and the relations \begin{enumerate}
\item $\langle a \rangle = \langle a b^2 \rangle$ for all $a,b \in k^\times$.
\item\label{atimesb} $\langle a \rangle \langle b \rangle = \langle a b \rangle$ for all $a,b \in k^\times/(k^\times)^2$
\item\label{aplusb} $\langle a \rangle + \langle b \rangle = \langle a+b \rangle + \langle a b (a+b) \rangle$ for all $a,b \in k^\times$ such that $a+b \neq 0$.
\end{enumerate} The relations \eqref{atimesb} and \eqref{aplusb} imply that
$$\langle a \rangle + \langle -a \rangle = \langle 1 \rangle + \langle -1 \rangle =: h$$ for all
$a \in k^\times/(k^\times)^2$, where $h$ is the hyperbolic form, defined $h = \langle 1 \rangle + \langle -1 \rangle$. 

As mentioned above, the generator $\langle a \rangle$ corresponds to the bilinear form $k \times k \to k$ taking $(x,y)$ to $a xy$. The fact
that $\{\langle a \rangle: a \in k^\times/(k^\times)^2 \}$ is a set of generators is thus equivalent to the statement from
bilinear algebra that a symmetric bilinear form can be stably diagonalized, meaning that after potentially adding
another form (which is only necessary in characteristic $2$), the original form is a perpendicular direct sum of bilinear forms on
$1$-dimensional vector spaces.

The {\em rank} of the bilinear form $\beta: V \times V \to k$ is the dimension of $V$ as a $k$-vector space. So for
example, the generators $\langle a \rangle$ for $a$ in $k^\times/(k^\times)^2$ are the isomorphism classes of all non-degenerate
rank one forms. The rank defines a ring homomorphism $$\rk: \GW(k) \to \bbb{Z}.$$ The kernel of the rank
homomorphism $$I:=\Ker (\rk)$$ is called the fundamental ideal and gives rise to a filtration of $\GW$
$$\GW \supseteq I \supseteq I^2 \supseteq \ldots$$ related to \'etale cohomology and Milnor K-theory by the Milnor
Conjecture \cite{MilnorAlgebraicKtheoryquadraticforms1970}, which is now a theorem of Voevodsky proven using
$\bbb{A}^1$-homotopy theory \cite{VoevodskyMotiviccohomologymathbfZ2003} \cite{VoevodskyReducedpoweroperations2003}.
Specifically, assume that the characteristic of $k$ is not $2$ and let $\rH^n(k, \bbb{Z}/2)$ denote the \'etale
cohomology of $\Spec k$ with $\bbb{Z}/2$ coefficients, or equivalently, the (continuous) cohomology of the Galois group
of the separable closure of $k$ with $\bbb{Z}/2$ coefficients. Let $ \KMf_n(k)$ denote the degree-$n$ summand of the
Milnor K-theory of $k$, \cite{MilnorAlgebraicKtheoryquadraticforms1970}. For present purposes, we only state that $\KMf_n(k)$
is the quotient of the $n$-fold tensor product $k^\times \otimes k^\times \otimes \ldots \otimes k^\times$ by the group
generated by tensors $a_1 \otimes a_2 \otimes \ldots \otimes a_n$ where there is some $i$ such that $a_i + a_{i+1} =
1$.
In degrees $0$, $1$ one has $\KMf_0(k) = \ZZ$ and $\KMf_1(k) = k^\times$. The Kummer map
$k^\times \to H^1(k, \bbb{Z}/2)$ is the first boundary map obtained by taking the Galois cohomology of the exact
sequence
$$ 1 \to \{ \pm 1\} \to k^\times_s \stackrel{z \mapsto z^2}{\longrightarrow} k^\times_s \to 1,$$ where $k_s$ denotes the separable
closure of $k$. It extends to a ring homomorphism $\KMf_*(k) \otimes \bbb{Z}/2 \to \rH^*(k, \bbb{Z}/2)$. The Milnor
conjecture says there are isomorphisms $$I^n/I^{n+1} \cong \KMf_n(k) \otimes \bbb{Z}/2 \cong \rH^n(k, \bbb{Z}/2),$$ where
the second isomorphism is induced by the map just discussed coming from the Kummer map and the first sends
$$a_1 \otimes a_2 \otimes \ldots \otimes a_n $$ in $\KMf_n(k)$ to
$$(\langle a_1 \rangle - \langle 1 \rangle) (\langle a_2 \rangle - \langle 1 \rangle) \cdots (\langle a_n \rangle -
\langle 1 \rangle)$$ in $I^n/I^{n+1}$.

The maps $I^n \to I^n/I^{n+1} \to \rH^n(k, \bbb{Z}/2)$ define invariants on bilinear forms in $I^n$. The $0$th of these
is by definition the rank. The first is the {\em discriminant}
\[\disc : \GW(k) \to \rH^1(k, \bbb{Z}/2) \cong k^\times/(k^\times)^2,\] which sends the bilinear form $\beta: V \times V \to k$ to
the determinant of a Gram matrix for $\beta$ $$\disc(\beta)=(\beta(v_i,v_j))_{i,j}$$ where $\{v_1, \ldots v_n\}$ is a
basis of $V$. The second is the Hasse--Witt invariant
$$\gamma: \GW(k) \to \rH^2(k, \bbb{Z}/2) \cong \KMf_2(k) \otimes \bbb{Z}/2,$$ which takes
$\langle a_1 \rangle + \langle a_2 \rangle + \ldots +\langle a_m \rangle$ to the image of $$\sum_{i<j} a_i \otimes a_j$$
in $\KMf_2(k) \otimes \bbb{Z}/2$. If we identify $\rH^2(k, \bbb{Z}/2) $ with the $2$-torsion subgroup of the Brauer group
of equivalence classes of central simple $k$-algebras, the class
$\gamma(\langle a_1 \rangle + \langle a_2 \rangle + \ldots +\langle a_m \rangle)$ is that of the tensor product of the
quaternion algebras $(a_i,a_j)$ as $(i,j)$ runs over pairs of integers from $1$ to $m$ with $i<j$. The third is the Arason invariant. The higher
invariants are not explicitly named, but do not vanish in general. For example, in the Grothendieck--Witt group of
the field $\bbb{C}((x_1,x_2,\ldots, x_n))$, the $n$th quotient $I^n/I^{n+1}$ is not $0$.

For many fields $k$, there are explicit computations of $\GW(k)$ and algorithms for deciding when two sums of the given generators represent the same element. See \cite{LamIntroductionquadraticforms2005}.

\begin{example}
When $k = \bbb{C}$ or any algebraically closed field, $k^\times/(k^\times)^2$ is the one element group, so by the previously given presentation of $\GW(k)$, we see that $\rk: \GW(k) \to \bbb{Z}$ is an isomorphism, e.g., $\GW(\bbb{C}) \cong \bbb{Z}$.
\end{example}

\begin{example}
 For a bilinear form $\beta: V \times V \to \bbb{R}$ over the real numbers, Sylvester's law of inertia states that
 there is a basis of $V$ so that $\beta$ is diagonal with only $1$'s, $-1$'s and $0$'s on the diagonal. The {\em
 signature} of $\beta$ is the number of $1$'s minus the number of $-1$'s. The signature determines a homomorphism
 $\sgn: \GW(k) \to \bbb{Z}$. Non-degenerate bilinear forms over $\bbb{R}$ are classified by their rank and signature,
 and thus the homomorphism $$\rk \times \sgn: \GW(\bbb{R}) \to \bbb{Z} \times \bbb{Z}$$ determines an isomorphism from
 $\GW(\bbb{R})$ onto its image $\GW(\bbb{R}) \cong \bbb{Z} \times \bbb{Z}$.
\end{example}

\begin{example}
For a finite field $\bbb{F}_q$, the rank and the discriminant define a group isomorphism $\GW(\bbb{F}_q) \cong \bbb{Z} \times \bbb{F}_q^*/(\bbb{F}_q^*)^2$. When $q$ is odd, note that this means that the Grothendieck--Witt group is isomorphic to $ \bbb{Z} \times \bbb{Z} /2$.
\end{example}

Springer's theorem \cite[VI Theorem 1.4]{LamIntroductionquadraticforms2005} computes the Grothendieck--Witt group of $\bbb{Q}_p$ and $\bbb{C}((x_1,x_2,\ldots,x_n))$. A theorem of Milnor gives the analogous computation for function fields \cite[Theorem 5.3]{MilnorAlgebraicKtheoryquadraticforms1970}. For $\GW(\bbb{Q})$, see \cite[VI 4]{LamIntroductionquadraticforms2005}. In summary, the Grothendieck--Witt group is interesting and computable.

Morel's degree has the agreeable feature that for subfields $k$ of $\bbb{R}$, it records the ordinary topological degrees of both the complex and real realizations. Namely, the following diagram commutes. $$\xymatrix{ \Top \ar[d]^{\deg} &\ar[d]^{\deg} \Spaces \ar[l]_{|\cdot|_{\RR}} \ar[r]^{|\cdot|} & \ar[d]^{\deg} \Top \\
 \bbb{Z} & \ar[l]_{\rk} \GW(k) \ar[r]^{\sgn}& \bbb{Z}} $$ 
 
 \subsection{Local Degree}
 
As in the case of the degree morphism from classical algebraic topology, there is an associated notion of ``local degree'',
defined as follows. Let $U$ and $W$ be Zariski open subsets of $\bbb{A}^n_k$, and let $f: U \to W$ be a map. Suppose $x$ is a closed point of $U$ such that $y=f(x)$
is a $k$-point (meaning that the coordinates of $y$ are elements of $k$), and such that $x$ is isolated in
$f^{-1}(y)$. By possibly shrinking $U$, we may assume that $f$ maps $U -\{x\}$ to $W - \{y\}$, and therefore $f$ induces
a map $$f: U/(U-\{x\}) \to W/(W-\{y\}).$$ The Purity Theorem identifies the quotients $W/(W - \{y\})$ and $U/(U-\{x\})$ with the
sphere $(\bbb{P}^1)^{\wedge n}$ and the smash product $(\bbb{P}^1)^{\wedge n}\wedge \Spec k(x)_+$, respectively. The
{\em local degree} $\deg_x f$ of $f$ at $x$ is the degree of the composite
$$ \bbb{P}^n/\bbb{P}^{n-1} \to \bbb{P}^n/(\bbb{P}^n -\{ x\}) \cong U/(U-\{x\}) \stackrel{f}{\longrightarrow} W/(W-\{y\})
\cong (\bbb{P}^1)^{\wedge n} \cong \bbb{P}^n/\bbb{P}^{n-1} $$
in $\GW(k)$. If $x$ is also a $k$-point, it is more natural to take the degree of the composite
$$(\bbb{P}^1)^{\wedge n} \to U/(U-\{x\}) \stackrel{f}{\longrightarrow} W/(W-\{y\}) \to(\bbb{P}^1)^{\wedge n} .$$ These
definitions are equivalent \cite[Proposition 11]{KWA1degree}. In the presence of Nisnevich local coordinates (see
\cite[Definition 17]{CubicSurface}) and appropriate orientations, the definition of local degree can be
generalized.

There is an explicit algorithm to compute the local degree \cite{eisenbud77} \cite{eisenbud78} \cite{KWA1degree}, which
can be implemented with a computer algebra package, or even by hand in good circumstances. Namely, consider a map
$f: U \to W$, and a point $x$ in $U$ as above. Since $W$ is a subset of $\bbb{A}^n_k$, we can express $f$ as an
$n$-tuple of functions $f=(f_1,\ldots, f_n)$ defined on $U$. Since $U$ is an open subset of $\bbb{A}^n_k$, the point $x$
corresponds to a prime ideal $p$ in $k[x_1,\ldots, x_n]$. The $f_i$ can be viewed as elements of the localization
$k[x_1,\ldots,x_n]_p$. The assumption that $x$ is isolated in $f^{-1}(y)$ implies that the quotient
$$Q = \frac{k[x_1,\ldots,x_n]_p}{\langle f_1,\ldots, f_n \rangle} $$ is finite dimensional as a $k$-vector space. The $\bbb{A}^1$-local degree of $f$ at $p$ is represented by the following bilinear form on $Q$. The
Jacobian $$J = \det \Big(\frac{\partial f_i}{\partial x_j}\Big)_{i,j} $$ determines an element of $Q$. For simplicity, we assume that $p$ is a $k$-point, and that $\rk_k Q$ is not divisible by the characteristic of $k$. Choose any $k$-linear function $$\eta: Q \to k $$ such that
$$\eta(J) = \rk_k Q.$$ Then define the bilinear form $\omega^{\EKL}: Q \times Q \to k$ by
$$\omega^{\EKL}(g,h) = \eta(gh).$$ It can be shown that the isomorphism class of $\omega^{\EKL}$ does not depend on the
choice of $\eta$ \cite[Proposition 3.5]{eisenbud77} and that the class of $\omega^{\EKL}$ in $\GW(k)$ is the $\bbb{A}^1$-local degree \cite{KWA1degree}. For the analogous algorithm without the two assumptions, see the work of Scheja and Storch \cite[Section 3]{scheja}. 

As in classical algebraic topology, the degree of a map can be expressed as a sum of local degrees. For example, let
$f: \bbb{P}^n_k \to \bbb{P}^n_k$ be a finite map. Let $\bbb{A}^n_k$ denote the Zariski open subset of $\bbb{P}^n_k$
which is the complement of the standard closed immersion $\bbb{P}^{n-1}_k \hookrightarrow \bbb{P}^n_k $ described in
Example \ref{Pn/Pn-1=Sn+nalpha}. Suppose that $f^{-1}(\bbb{A}^n_k) = \bbb{A}^n_k$, so there is an induced map
$$F: \bbb{P}^n_k/\bbb{P}^{n-1}_k \to \bbb{P}^n_k/\bbb{P}^{n-1}_k.$$ Then, for any $k$-point $y$ of
$\bbb{A}^n_k$, \begin{equation*} \deg F = \sum_{x \in f^{-1}(y)} \deg_x f .\end{equation*} This is well-known and the
proof is very similar to a classical argument. A reference is \cite[Proposition 13]{KWA1degree}.
 
\subsection{Milnor--Witt groups}
The ring $\GW(k)$ is the degree $0$ subring of a graded ring $\KMW_{\ast}(k)$ developed by Morel and M.~J.~ Hopkins \cite{MorelmathbbAalgebraictopology2012} to describe part of the motivic homotopy groups of
spheres. For a field $k$, the ring $\KMW_{\ast}(k)$ is the associative, graded ring with the following presentation. The
generators consist of one element $\eta$ of degree $-1$ and elements $[a]$ of degree $1$ for each $a$ in $k^\times$. The
hyperbolic element $h$ is defined $h:=\eta [-1] + 2$. These generators are subject to the relations \begin{itemize}
\item (Steinberg relation) $[a][1-a] = 0$ for each $a$ in $k -\{0,1\}$. 
\item $[ab] = [a] + [b] + \eta [a][b]$ for every $a,b$ in $k^\times$.
\item $[a] \eta = \eta [a]$ for every $a$ in $k^\times$.
\item $\eta h = 0$.
\end{itemize} One can identify $\GW(k)$ with $\KMW_{0}(k)$ by identifying $\langle a \rangle$ with $1+ \eta [a]$,
$$ \langle a \rangle = 1+ \eta [a],$$ and it follows that the two definitions $h=\eta [-1] + 2$ and
$h = \langle 1 \rangle + \langle -1 \rangle$ agree. Morel proves an isomorphism between $\bbb{P}^1$-stable
$\bbb{A}^1$-homotopy classes of maps and $\KMW_{\ast}(k)$ $$[S^0, \bbb{G}_m^{\wedge n}]_{\soneSpt} \cong \KMW_{n}(k).$$
Moreover, there is an analogous unstable result:
$$ [S^n \wedge \bbb{G}_m^{\wedge r}, S^n \wedge \bbb{G}_m^{\wedge j}] \cong \KMW_{j-r} (k)$$ for $n \geq 2$, which
extends to an isomorphism of sheaves discussed in Section \ref{section:connectivity_theorem}. For $n=1$, we are lead to
consider the unstable $\bbb{A}^1$-homotopy classes of maps $\bbb{P}^1_k \to \bbb{P}^1_k$, also discussed in Section
\ref{section:connectivity_theorem}; see Equation~\eqref{unstableA1homotopyclassesP1toP1}.

The fact that $\GW(k)$ and $\KMW_{\ast}(k)$ are global sections of stable $\Aone$-homotopy sheaves
$\bpi_{\ast}^{\bbb{A}^1\!,s}$ implies the existence of extra structure, for example
residue maps, which are closely related to the unramified sheaves discussed in more detail in the next section, and transfer maps, which exist for $2$-fold $\G_m$-loop
sheaves of $\bbb{A}^1$-homotopy sheaves $\pi_{\ast}^{\bbb{A}^1}$ for $\ast \geq 2$, and in particular exist for stable
$\bbb{A}^1$-homotopy sheaves \cite[Chapter 4]{MorelmathbbAalgebraictopology2012}). We discuss these residue and transfer
maps in more detail now.

We have a residue map in the following situation. Let $v: K \to \bbb{Z} \cup \infty$ be a valuation on a field $K$
containing $k$, and let $\calO_v$ denote the corresponding valuation ring of elements of non-negative valuation. The
residue field $k(v)$ is the quotient of $\calO_v$ by its maximal ideal consisting of the elements of $K$ of positive
valuation. Let $\pi$ be a uniformizer (i.e., an element such that $v(\pi)=1$). Then there is a unique homomorphism
$\partial_v^{\pi}$ called a {\em residue map}\begin{equation*} \partial_v^{\pi}: \KMW_{\ast}(K) \to \KMW_{\ast -1}(k(v))
\end{equation*} commuting with $\eta$ and such that $$\partial_v^{\pi} ([\pi][a_1][a_2]\ldots[a_n]) =
[\overline{a_1}][\overline{a_2}]\ldots[\overline{a_n}],$$ $$\partial_v^{\pi} ([a_0][a_1][a_2]\ldots[a_n]) = 0 ,$$ where
$a_i\in \calO_v^*$ and $\overline{a_i}$ denotes the image of $a_i$ in $k(v)^*$. The kernel of the residue maps define
the sections of the sheaf $\KMW_{\ast}$ over $\Spec \calO_v$, see the discussion following
\cite{MorelmathbbAalgebraictopology2012}[Lemma 3.19].

We have a transfer map in the following situation. Let $K \subseteq L$ be a field extension of finite rank, where $K$ is finite type over $k$. 

We first note that the
inclusion $K \subseteq L$ corresponds to a map of schemes (or spaces) $\Spec L \to \Spec K$ in the opposite
direction. The sheaf property (which comes from pullback of maps) gives restriction maps
$\KMW_{\ast}(K) \to \KMW_{\ast}(L)$. These restriction maps are, not surprisingly, $[a] \mapsto [a]$,
$\eta \mapsto \eta$, and correspond to pullback of bilinear forms when restricted to $\GW$.

Transfer maps for such field extensions can be constructed by producing a stable map $$\Spec K \to \Spec L$$ in the
direction opposite to the map of spaces, closely analogous to the Becker--Gottlieb
transfer\cite[4.2]{MorelmathbbAalgebraictopology2012}. Namely, when $L$ is generated by a single element over $K$, we
can choose a closed point $z$ of $\bbb{P}^1_K$ with residue field $L$, or equivalently, the data of a closed
immersion $$z: \Spec L \hookrightarrow \bbb{P}^1_K.$$ Using Purity, we then have the cofiber sequence
\[(\bbb{P}^1_K-\{z\}) \to \bbb{P}^1_K \to \bbb{P}^1_K/(\bbb{P}^1_K-\{z\}) \cong \Th (N_z \bbb{P}^1_K) \cong \bbb{P}^1_k
 \wedge (\Spec L)_+.\]

The quotient map, or Thom collapse map, in this sequence
\[\mathbb{P}_K^1 \cong \bbb{P}^1_k \wedge (\Spec K)_+ \to \bbb{P}^1_k \wedge (\Spec L)_+\] is the desired stable map
$\Spec K \to \Spec L$, and induces a transfer map $$\tau_{L/K}^z:\KMW_{\ast}(L) \to \KMW_{\ast}(K),$$ called the {\em
 geometric transfer}. This transfer depends on the chosen generator of $L$ over $K$. Although $L$ may not be generated
over $K$ by a single element, we can always choose a finite list of generators. Given an ordered such list, we define a
transfer $\KMW_{\ast}(L) \to \KMW_{\ast}(K)$ by composing the transfers just constructed.
 
For simplicity, assume that the characteristic of $k$ is not $2$. The dependency of the transfer map on the chosen
generators can be eliminated by modifying the definition in the following manner: Suppose $z$ is a generator of $L$ over
$K$. The monic minimal polynomial of $z$ can be expressed in a canonical manner as $P(x^{p^m})$ where $P$ is a separable
polynomial. The derivative $P'(z^{p^m})$ of $P$ evaluated at $z^{p^m}$ is an element of $L^*$. The {\em cohomological
 transfer} \begin{equation}\label{cohtrmap} \Tr_{L/K}: \KMW_{\ast} (L) \to \KMW_{\ast}(K)\end{equation} is then
defined $$\Tr_{L/K}(\beta) = \tau_{L/K}^z (\langle P'(z^{p^m})\rangle \beta),$$ and this map is independent of the
chosen generator \cite[Theorem 4.27]{MorelmathbbAalgebraictopology2012}. More generally, for any finite extension
$K \subseteq L$, one has a cohomological transfer \eqref{cohtrmap} by composing the the transfers just defined for a
sequence of generators, and the resulting map is again independent of the chosen sequence. This independence of the
choice of generators is more naturally understood terms of twisted Milnor--Witt K-theory, where the twist is by the
canonical sheaf.

From the algebraic perspective, there are many possible transfers $\GW(L) \to \GW(K)$ for $[L : K] \leq \infty$ as
follows: for any nonzero $K$-linear map $f:L \to K$, and non-degenerate bilinear form $\beta: V \times V \to L$, the
composite $$f \circ \beta:V \times V \stackrel{\beta}{\to} L \stackrel{f}{\to} K $$ is a non-degenerate bilinear
form. Thus we may define a map $\Tr_f: \GW(L) \to \GW(K)$ which takes the isomorphism class of an appropriate $\beta$ to
the isomorphism class of $f \circ \beta$, where in the former case, $V$ is viewed as an $L$-vector space, and in the
latter case, $V$ is viewed as a $K$-vector space.

This abundance of transfers for $\GW$ also implies the same for $\KMW_{\ast}$ as the latter can be expressed as the
fiber product $$\xymatrix{\KMW_{\ast}(K) \ar[r] \ar[d] &I(K)^n \ar[d]\\ \KMf_n(K) \ar[r] & I(K)^{n}/I(K)^{n+1},}$$ and
there is a canonical transfer on Milnor K-theory.

When $K \subset L$ is a finite separable extension, there is a canonical choice of such an $f$, namely, the trace map
$L \to K$ from Galois theory, given by summing the Galois conjugates of an element of $L$. The resulting transfer map is
the cohomological transfer \eqref{cohtrmap} defined above. This transfer arises naturally when studying the local degree
at non-$k$-rational points. For example, let $f: U \to W$ be a map between open subsets of $\bbb{A}^n_k$. Suppose that
$x$ in $U$ maps to a $k$-point $y$ of $W$ and that $x$ is an isolated in $f^{-1}(y)$, so the local degree
$\deg^{\mathbb{A}^1}_x f$ exists. If the extension $k \subseteq k(x)$ is separable and $f$ is \'etale at $x$, then
$$ \deg^{\mathbb{A}^1}_x f = \Tr_{k(x)/k} \langle J(x)\rangle,$$ where $J$ denotes the Jacobian of $f$. A proof of this
is in \cite[Proposition 14]{KWA1degree}, and this proof relies on Hoyois's work in \cite{Hoyois_lef}.

\begin{example}
$\Tr_{\bbb{C}/\bbb{R}} \langle 1 \rangle = \begin{bmatrix} 2 &0 \\0 & -2\end{bmatrix} = \langle 1 \rangle + \langle -1 \rangle= h,$ where the central expression is the Gram matrix with respect to the $\bbb{R}$-basis $\{1,i\}$ of $\bbb{C}$.
\end{example}

\subsection{Euler class}\label{EulerClass_section}

Given an oriented vector bundle of rank $r$ on an oriented $\mathbb{R}$-manifold of dimension $r$, one can define an Euler number. It is an element of $\bbb{Z}$ that counts the number of zeros of a section in the following sense: Given a section $\sigma$ with an isolated zero $x$, we can define a local index (or degree) of $\sigma$ at $x$ by the following procedure. First, choose local coordinates around $x$ and a local trivialization
of the vector bundle such that both are compatible with the appropriate orientations. With these choices, the section
can be viewed as a function $\mathbb{R}^r \to \mathbb{R}^r$. The index $\ind_x \sigma$ is then the local degree of this
function at $x$, where $x$ is viewed as a point of $\mathbb{R}^r$, by a slight abuse of notation. The Euler number $e$ is then \begin{equation}\label{e=sumindex} e= \sum_{x: \sigma(x) = 0} \ind_x \sigma,\end{equation} when $\sigma$ has only isolated zeros.

With Morel's $\GW(k)$-valued degree, this Euler number can be enriched to an element of $\GW(k)$. Namely, one can again
define a local index using the local degree and then use \eqref{e=sumindex} to define the Euler number, now an element
of $\GW(k)$ \cite{CubicSurface}. There are subtleties involved in choosing coordinates and identifying $\sigma$ with a
function, due to difficulties in defining algebraic functions. One must then show that the Euler number is
well-defined. There are other approaches avoiding these difficulties, using oriented Chow or Chow--Witt groups of Barge
and Morel \cite{bargeGroupeChowCycles2000}. Barge and Morel construct an Euler class (loc. cit.) in oriented Chow and it
can be pushed forward under certain conditions to land in $\GW(k)$. See the work of Jean Fasel
\cite{faselGroupesChowWitt2008} and Marc Levine \cite{Levine-EC}. Morel gives an alternate construction of an Euler
class using obstruction theory \cite[Section 8.2]{MorelmathbbAalgebraictopology2012}. This Euler class lies in the same
oriented Chow group, and will be discussed further in Section \ref{Section:splittingproblem}. Some comparison results
are available between these two Euler classes \cite{asokComparingEulerClasses2016} \cite[Proposition
11.6]{Levine-EC}. Further approaches to defining an Euler class using $\bbb{A}^1$-homotopy theory are found in
\cite{DJK} \cite{LevineRaksit_MotivicGaussBonnet}.

\subsection{Applications to enumerative geometry} \label{Applications_to_enumerative_geometry_section}

Questions in enumerative geometry ask to count a set of algebro-geometric objects satisfying certain conditions, or more
generally, to describe this set. A classical example is the question ``How many lines intersect four general lines in
three dimensional space?" The questions are typically posed so there is a fixed answer, such as ``$2$," as opposed to
``sometimes $2$ and sometimes $0$," and to get such ``invariance of number," one needs to work over an algebraically
closed field; the number of solutions to a polynomial equation of degree $n$ is always $n$ if one works over $\bbb{C}$,
but not over $\bbb{R}$, and the same phenomenon appears in the classical question quoted above about the number of
lines. However, a feature of $\bbb{A}^1$-homotopy theory is its applicability to general fields $k$. Moreover, there are
classical theorems in enumerative geometry where a count of geometric objects is identified with an Euler number, which
is then computed using characteristic class techniques. Equipped with an enriched Euler class in $\GW(k)$, as in 
Section \ref{EulerClass_section}, we may thus take these theorems and hope to perform the analogous counts over other
fields.

What sorts of results appear? Here is an example taken from \cite{CubicSurface}. A cubic surface $X$ is the zero locus
in $\bbb{P}^3_k$ of a degree $3$ homogenous polynomial $f$ in $k[x_0,x_1,x_2,x_3]$. It is a lovely $19$th century result
of Salmon and Cayley that when $k = \bbb{C}$ and $X$ is smooth, the number of lines in $X$ is exactly $27$. Over the
real numbers, the number of real lines is either $3$, $7$, $15$, or $27$, and in particular, the number depends on the
surface. A classification was obtained by Segre \cite{segre42}, but it is a recent observation of Benedetti--Silhol
\cite{Benedetti}, Finashin--Kharlamov \cite{FinashinAbundanceRealLines2013}, Horev--Solomon \cite{solomon12} and Okonek--Teleman
\cite{okonek14} that a certain signed count of lines is always $3$. Specifically, Segre distinguished between two types
of real lines on $X = \{ [x_0,x_1,x_2,x_3] \in \bbb{R}\bbb{P}^3: f(x_0,x_1,x_2,x_3) = 0 \}$, called {\em hyperbolic} and
{\em elliptic}. The distinction is as follows. Let $L$ be a real line on $X$. For every point $p$ of $L$, consider the
intersection $T_p X \cap X$ of the tangent plane to the cubic surface at $p$ with the cubic surface itself. Since
$T_p X$ is a plane, the intersection is a curve in the plane, which by B\'Bezout's theorem has degree $3 \cdot 1
=3$.
The line $L$ must be contained in this plane curve. Algebraic curves can be decomposed into irreducible components, and
it follows that we can express $T_p X \cap X$ as a union $T_p X \cap X = L \cup Q$, where $Q$ is a plane curve of degree
$3 - 1 = 2$. Applying B\'ezout's theorem again shows that the intersection $L \cap Q$ of the plane curves $L$ and $Q$ is
degree $2$, or in other words, generically consists of two points, which we call $\{p,q\}$. See Figure
\ref{Cubic_surface_involution_picture}. We may thus define an involution $I: L \to L$ by sending $p$ to the unique point
$q$.

\begin{figure}
\begin{tikzpicture}
 \draw[thick] (-2, 0) -- (2, 0) node[right] {$L$};
 \draw[thick] (0, 0) circle[radius=1];
 \node[above right] at (.6, .6) {$Q$};
 \node[point] at (-1, 0) {};
 \node[below right] at (1, 0) {$q$};
 \node[point] at (1, 0) {};
 \node[below left] at (-1, 0) {$p$};

\end{tikzpicture}
\begin{caption}\relax
$T_p X \cap X \subset T_p X$ 
\label{Cubic_surface_involution_picture}
\end{caption}
\end{figure}

The points of the intersection $L \cap Q=\{p,q\}$ are precisely the points $x$ on $L$ such that $T_{x} X = T_{p} X$. To see this, note that $T_q X$ contains the span of a vector along $L$ and a vector along $Q$. At least generically, it follows that $T_q X$ contains a $2$-dimensional subspace of $T_p X$. Since $T_q X$ is $2$-dimensional ($X$ is smooth), it follows that $T_q X = T_p X$. Similar reasoning applies in the reverse direction as well. So we can characterize the involution $I$ as the unique map exchanging points on $L$ with the same tangent space to $X$.

An automorphism of $L \cong \mathbb{R}\bbb{P}^1$ is a conjugacy class of element of $\bbb{P}\GL_2 \bbb{R}$, and the elements of $\bbb{P}\GL_2 \bbb{R}$ are classified as elliptic, hyperbolic, or parabolic by the behavior of the fixed points. If $I$ has a complex conjugate pair of fixed points, $I$ is elliptic, if $I$ has two real fixed points, it is hyperbolic, and if $I$ has one fixed point, it is parabolic. Involutions are never parabolic. Segre classified the line $L$ as elliptic (respectively hyperbolic) if $I$ is. 

Another description of this distinction involves (S)pin structures. The tangent plane $T_p X$ rotates around $L$ as $p$ travels along $L$, describing a loop in the frame bundle of $\bbb{P}^3$. $L$ is hyperbolic if this loop lifts to the double cover, and hyperbolic if it does not.

The signed count referred to above is that \begin{equation}\label{hyperbolic-elliptic=3} \# \{ \text{hyperbolic lines} \}- \# \{ \text{elliptic lines} \}= 3.\end{equation}

Results which are true over $\bbb{C}$ and $\bbb{R}$ may be realizations of a more general result in $\bbb{A}^1$-homotopy theory. In the case of the count of lines on cubic surface, this is indeed the case: Consider a cubic surface $X$ over a field $k$. Suppose $L$ is a line in $\bbb{P}^3_{\overline{k}}$ which lies in $X$. The coefficients of $L$ determine a field extension $k(L)$. Moreover, the previously given definition of the involution $I$ carries over in this generalized situation, determining an involution $I$ of the line, and thus a conjugacy class in $\bbb{P}\GL_2 k(L)$. Such a conjugacy class has a well-defined determinant $\det I$ in $k(L)^*/(k(L)^*)^2$. The {\em Type} of $L$ is $\Type (L) = \langle \det I \rangle$ in $\GW(k(L))$. Alternatively, the type of $L$ may be described as $\Type (L) = \langle D \rangle$, where $D$ is the unique element of $k(L)^*/(k(L)^*)^2$ so that the fixed points of the involution $I$ are a conjugate pair of points defined over the field $k(L)[\sqrt{D}]$. There is a third description of the type as $\Type (L) = \langle -1 \rangle \deg I$ the multiplication of the $\bbb{A}^1$-degree of $I$ and $\langle -1 \rangle$ in $\GW(k(L))$. The theorem of Salmon and Cayley and \eqref{hyperbolic-elliptic=3} then are realizations of the following theorem \cite[Theorem 1]{CubicSurface}.

\begin{theorem}
Let $X$ be a smooth cubic surface over a field $k$ of characteristic not $2$. Then $$\sum_{\text{lines } L\text{ in } X} \Tr_{k(L)/k} \Type(L) = 15 \langle 1 \rangle + 12 \langle -1 \rangle.$$
\end{theorem}

This is proven by identifying the left hand side with the Euler class of the third symmetric power of the dual tautological bundle on the Grassmannian of lines in $\bbb{P}^3$.

Other results along these lines include \cite{Hoyois_lef}, \cite{Levine-EC}, \cite{Levine-Witt}, \cite{Levine-Welschinger}, \cite{Wendt-oriented_schubert}, \cite{SrinivasanArithmeticCountLines2018}, \cite{Bethea}, and this is an active area of research.

\section{Homotopy sheaves and the connectivity theorem}\label{section:connectivity_theorem}

In the monograph \cite{MorelmathbbAalgebraictopology2012}, Morel establishes a number of extremely strong results
describing the $\Aone$-homotopy sheaves of objects of $\Spaces$ when $k$ is an infinite perfect field. These are
sheaves on the big Nisnevich site $\Sm_k$. The assumptions on the field are probably both unnecessary, but the
literature does not currently contain proofs of certain necessary statements for finite or imperfect fields.

The $\Aone$-homotopy sheaves $\bpi_i^{\Aone}(X)$---when $i \ge 1$---must satisfy a strong $\Aone$-invariance
property. For all smooth $k$-schemes $U$, the map on cohomology
\begin{equation}
 \label{eq:2}
 \Hoh^n ( U , \bpi_i^{\Aone}(X) ) \to \Hoh^n ( U \times \Aone_k, \bpi_i^{\Aone}(X) ) 
\end{equation}
must be an isomorphism for all applicable $n$, to wit, $n\in \{0,1\}$ for all sheaves of groups, and
$n \in \{ 0, 1 \dots, \}$ if the group is abelian. This condition is known as \textit{strong $\Aone$-invariance} in the
general case, and \textit{strict $\Aone$-invariance} in the case of an abelian sheaf of groups. Over an infinite perfect
field, a strongly $\Aone$-invariant sheaf of abelian groups is strictly $\Aone$-invariant, \cite[Theorem
5.46]{MorelmathbbAalgebraictopology2012}, so for us the distinction is chiefly useful for defining a category of
``strictly $\Aone$-invariant sheaves'', which necessarily consist of abelian groups.

We remark that the case of $n=0$ in \eqref{eq:2} is merely $\Aone$-invariance. It is possible to give an example of an
(abelian) group that is $\Aone$-invariant but not strongly $\Aone$-invariant---for instance, if one defines
$\tilde\ZZ[\Gm]$ to be the abelian group generated by the sheaf of sets $\Gm$, modulo the relation that $1 \in \Gm$ is
identified with the $0$ element, then $\tilde \ZZ [\Gm]$ is $\Aone$-invariant---since $\Gm$ itself is $\Aone$
invariant---but the free strongly $\Aone$-invariant sheaf of abelian groups generated by $\tilde \ZZ[\Gm]$ is $\KMW_1$
by \cite[Theorem 3.37]{MorelmathbbAalgebraictopology2012}.

A sheaf of groups, $\sh G$, is strongly $\Aone$-invariant if and only if it appears as an $\Aone$-homotopy sheaf
$\bpi^{\Aone}_1$, since it appears as the $\Aone$-fundamental sheaf of $B \sh G$. The strongly $\Aone$-invariant sheaves
$\sh G$ are \textit{unramified} sheaves, which is to say, briefly, that $\sh G(X)$ is the product of $\sh G(X_\alpha)$
as $X_\alpha$ range over the irreducible components, and that for any dense open $U \subset X$, $\sh G(X) \to \sh G (U)$
is an injection which is even an isomorphism if $X - U$ is everywhere of codimension at least $2$ in $X$, see
\cite[Section 2.1]{MorelmathbbAalgebraictopology2012}. A feature of unramified sheaves $\sh G$ in this sense is that $\sh G$ can be
recovered from the values $\sh G(F)$ where $F$ ranges over fields of finite transcendence degree over $k$, along with
subsets $\sh G (\sh O_v ) \subset \sh G (F)$ associated to any discrete valuation $v$ on $F$, and specialization maps
$\sh G (\sh O_v) \to \sh G( \kappa(v))$ mapping to the residue fields, all satisfying certain compatibility axioms,
\cite[Section 2.1]{MorelmathbbAalgebraictopology2012}. Most strikingly, if $\sh G(F)$ vanishes for all field extensions
of $k$, then $\sh G = 0$.

The category of strictly $\Aone$-invariant sheaves forms a full abelian subcategory of the category of sheaves of
abelian groups on $\Sm_k$, being the heart of a $t$-structure---this appears as \cite[Corollary
6.24]{MorelmathbbAalgebraictopology2012}. 

One can rephrase the statement that the sheaves $\bpi_i^{\Aone}(X)$ are strongly $\Aone$-invariant for $i \ge 1$ as
follows: if the $\Aone$-localization map $X \to L_{\Aone} X$ is a simplicial equivalence, i.e., if $X$ is $\Aone$-local,
then the sheaves $\bpi_i(X)$ for $i \ge 1$ are strongly $\Aone$-invariant. There is a partial converse:
\begin{theorem}
 Suppose $X$ is a pointed connected object of $\Spaces$. Then $X$ is $\Aone$-local if and only if $\bpi_i(X)$ is
 strongly $\Aone$-invariant for all $i \ge 1$.
\end{theorem}
This holds, loosely speaking, because $X$ is $\Aone$-local if the functor
$\Map(\cdot, X)$ does not distinguish between $U$ and $U \times \Aone$ up to homotopy. For a connected $X$, the calculation of maps
$\Map(\cdot, X)$ can be reduced to the calculation of $\Map(\cdot, K(\bpi_i(X), j))$ for varying $i, j$ by means of the
Postnikov tower, and strong $\Aone$-invariance of $\bpi_i(X)$ amounts to the same thing as $\Aone$-locality of
$K(\bpi_i(X), j)$. The result is stated as \cite[Corollary 6.3]{MorelmathbbAalgebraictopology2012}.

The assumption of connectivity, while it may seem mild at first, is a great inconvenience. In contrast to classical
homotopy theory, where one can generally work component-by-component, it is a feature of the homotopy of sheaves that
$\bpi_0$, a sheaf of sets, can be extremely complicated. 

In fact, the difficulties at $\bpi_0$ are substantial. If $X$ is a discrete sheaf of groups, viewed as a
space in dimension $0$, then $X$ is $\Aone$-local if and only if $X$ is $\Aone$-invariant. If $X$ is $\Aone$-invariant
without being strongly $\Aone$-invariant, then $L_{\Aone} B X$ cannot be (simplicially) equivalent to $B X$, and so
$ \bpi_0\Omega L_{\Aone} BX \not \iso \bpi_0 X$. In particular, $BX$ does not have the $\Aone$-homotopy type of a
delooping of $X$. It follows from this that $\Aone$-homotopy theory does correspond to an $\infty$--topos in the sense
of \cite[Section 6.1]{LurieHighertopostheory2009}. This observation, due to J.~Lurie, appears as \cite[Remark
3.5]{SpitzweckMotivictwistedtheory2012}.

As a heuristic, aside from problems at $\bpi_0$, results that hold in $\infty$-topoi (or in the ``model topoi'' of
C.~Rezk \cite{ToenHomotopicalalgebraicgeometry2005}) can generally be expected to hold in $\Aone$-homotopy theory. For
instance:
\begin{theorem}
 Suppose 
 \begin{equation}
 \label{eq:3}
 F \to E \to B
 \end{equation}
 is a simplicial homotopy fiber sequence in $\Spaces$ where $\bpi_0(B)$ is $\Aone$-invariant. Then \eqref{eq:3} is an
 $\Aone$-homotopy fiber sequence.
\end{theorem}
This is a corollary of \cite[Theorem 2.1.5]{AsokAffinerepresentabilityresults2018} (see also Remark 2.1.6), a
development of \cite[Theorem 6.53]{MorelmathbbAalgebraictopology2012}.

\subsection{Contractions}

The following construction appears originally in \cite{VoevodskyCohomologicaltheorypresheaves2000}, but applied there
only to
``presheaves with transfers''. In the current context, it is due to \cite{MorelmathbbAalgebraictopology2012}[Section
2.2, pp. 33--36]. Given a
presheaf of groups $\sh G$, define
\[ \sh G_{-1}(U) = \ker( \sh G(\Gm \times U) \overset{ev_1}{\longrightarrow} \sh G(U). \]
Here the map is induced inclusion of $U \times \{1\}$ in $U \times \Gm$; when the presheaf is applied, this appears as a
kind of ``evaluation at $1$'', hence the notation. This construction is functorial in $\sh G$.

The functor $\sh G \mapsto \sh G_{-1}$ is sometimes called \textit{contraction}. It may be iterated, in which case one
writes $\sh G_{-n}$. The functor $(\cdot)_{-1}$ has a number of excellent properties. It restricts to give a functor of
(Nisnevich) sheaves, for instance, and it preserves the property of being abelian. It is also left exact in
general. When applied to $\Aone$-homotopy sheaves, it has the following striking description, which appears as
\cite[Theorem 6.13]{MorelmathbbAalgebraictopology2012}
\begin{theorem}[Morel]
 If $X$ is a pointed, connected $\Aone$-local space, then the (derived) mapping space $\Omega_{\Gm} := \Map(\Gm, X)$ is also pointed,
 connected and $\Aone$-local, and there is a canonical isomorphism
 \[ \bpi_n^{\Aone}( \Omega_{\Gm} X) \to \bpi_n^{\Aone}(X)_{-1}. \]
\end{theorem}
This theorem, of course, is known only for homotopy sheaves over $\Sm_k$ where $k$ is a field satisfying the running
assumptions of \cite{MorelmathbbAalgebraictopology2012}.

Any short exact sequence of strictly $\Aone$-invariant sheaves may be realized as the homotopy sheaves of an $\Aone$
homotopy fiber sequence, for instance
\[0 \to \bpi_2^{\Aone}( K(\sh A,2)) \to \bpi^{\Aone}_2(K(\sh B, 2)) \to \bpi^{\Aone}_2(K(\sh C, 2)) \to 0 \]
and since the (derived) mapping space $\Map(\Gm, \cdot)$ preserves homotopy fibre sequences, one deduces that the
functor $\sh A \mapsto \sh A_{-1}$ is exact on the category of strictly $\Aone$-invariant sheaves.

It is immediate that if $\sh C$ is a constant presheaf of abelian groups, then $\sh C_{-1} = 0$. One might fantasize
that if $\sh G$ is strictly $\Aone$-invariant, then $\sh G_{-1} = 0$ should imply that $\sh G$ is constant, but this is
far from the case. In fact, from the analysis furnished by \cite[Chapter 2]{MorelmathbbAalgebraictopology2012} of
$(\cdot)_{-1}$, one can deduce that $\sh G_{-1} = 0$ if and only if $\sh G$ is \textit{birational} in the sense of
converting dense open inclusions of schemes $U \to X$ into isomorphisms $\sh G(X) \to \sh G(U)$. In \cite[Section
6]{AsokSmoothvarietiesBbb2011}, a study is made of such sheaves, and the category of all such sheaves over $\Sm_k$ is
seen to be equivalent to a very large category of functors on the category of field extensions of $k$. 

\subsection{Unramified $K$-theories}

Certain groups that were previously known to be functors of fields and field extensions are known to extend to give
strictly $\Aone$-invariant sheaves on the Nisnevich site of $\Sm_k$. For instance, there is a strictly $\Aone$-invariant
sheaf $\K^Q_n$ such that $\K^Q_n(F) = K^Q_n(F)$, Quillen's $K$-theory---here the field extension $F/k$ is supposed to be a
separable extension of an extension of finite transcendence degree. The sheaf $\K^Q_n$ arises as
$\bpi^{\Aone}_n(B \GL_N)$ for $N \ge n+2$, this being a consequence of the representability of $K$-theory for smooth
schemes as proved in \cite[Theorem 3.13, p140]{morelMathbbAHomotopyTheory1999}. Analogous constructions for hermitian
$K$-groups are made in \cite{SchlichtingGeometricmodelshigher2015}.

Another notable sheaf is $\K^M_n$, which can be recovered as a quotient of $\KMW_n$, the
\textit{unramified Milnor--Witt $K$-theory sheaf} as constructed by Morel in \cite[Section
3.2]{MorelmathbbAalgebraictopology2012}. In this case, the group of sections for a field is $K^M_n(F)$, as defined in
Section \ref{degree_section}.

In the cases of $\K^Q_n$, $\K^{MW}_n$ and $\K^M_n$, the phenomenon of $\PP^1$-stability for the associated theories
implies $(\K^Q_n))_{-1} = \K^Q_{n-1}$ and similarly for the other two theories; a proof is outlined in \cite[Lemma 2.7,
Proposition 2.9]{AsokAlgebraicvectorbundles2014}.

One may define $\K_3^{\text{ind}}$, the cokernel of a natural map $\phi: \K^M_3 \to \K^Q_3$. By virtue of Matsumoto's
theorem, one knows that
\[ (\phi)_{-1}: \KM_2 \to \K^Q_2 \]
is an isomorphism, so that $(\K^{\text{ind}}_3)_{-1} \iso 0$. It is known, for instance by
\cite{Merkurevgroupfield1990}, that $\K_3^{\text{ind}}$ is not constant, so this furnishes a specific example of a nonconstant
strictly $\Aone$-invariant sheaf, the contraction of which is $0$.

\subsection{$\Aone$-homology and the connectivity theorem}

In \cite{MorelmathbbAalgebraictopology2012}, Morel defines an ``$\Aone$-homology theory'' $\Hoh_n^{\Aone}(X)$, by means
of the following: for any simplicial (pre)sheaf $X$ one may define $\ZZ(X)$, the free abelian group on $X$, which one
converts to a chain complex $C_*(X)$ via the Dold--Kan correspondence. One can take the category of (pre)sheaves of
chain complexes on $\Sm_k$, viewed as a setting for homotopy theory in its own right, and then localize with respect to
$\Aone$, that is, with respect to the maps $C_* \tensor C_*(\Aone) \to C_*$. One may now replace chain complexes
$C_*(X)$ by ``abelian $\Aone$-local replacements''---for Morel, these are the fibrant objects in the localized model
category---denoted $C^{\Aone}_*(X)$. The homology of such an object is $\Hoh_*^{\Aone}(X)$, and these homology sheaves
are all strictly $\Aone$-invariant.

This homology theory is quite distinct from the ``motivic homology theory'' defined by Voevodsky \cite{VoevodskyCohomologicaltheorypresheaves2000}; it is much
closer to the unstable $\Aone$-homotopy theory, and little is known about it. It does enjoy the following three
properties:
\begin{enumerate}
\item It is $S^1$ stable: $\Hoh_{n+1}^{\Aone}( X \wedge S^1 ) \iso \Hoh_n^{\Aone}(X)$.
\item If $\sh F$ is a discrete sheaf, then $\Hoh^{\Aone}_0(\sh F)$ is the strictly $\Aone$-invariant sheaf freely
 generated by the sheaf $\sh F$, in the sense that this construction is left adjoint to an obvious forgetful functor.
\item For any $\Aone$-simply-connected pointed object $X$, a Hurewicz isomorphism holds. If $\bpi_i^{\Aone}(X) = 0$ for
 $i \le n-1$, and $n \ge 2$, then a natural map $\bpi_n^{\Aone}(X) \to \Hoh_n^{\Aone}(X)$ is an isomorphism. A
 modification of this holds for $\bpi_1$, involving abelianization and ``$\Aone$-strictification''. See \cite[Theorems
 6.35, 6.37]{MorelmathbbAalgebraictopology2012}.
\end{enumerate}

It is this theory, and the properties above, that allows Morel to compute the unstable $\Aone$-homotopy sheaves of the
spheres. Specifically, provided $n \ge 2$, so that $\bpi_n^{\Aone}$ is known to be abelian \& strictly $\Aone$
invariant, then $\bpi_n^{\Aone} (S^n \wedge \Gm^{\wedge m})$ is the free strictly $\Aone$-invariant sheaf generated by
the set $\Gm^{\wedge m}$. That is

\begin{equation}
 \label{eq:4}
 \bpi_n(S^n \wedge \Gm^{\wedge m}) = \KMW_m
\end{equation}
 provided $n \ge 2$.

Combined with further results of Morel's on the contractions of $\KMW_\bullet$, one deduces
\begin{equation}
 \label{eq:5}
 \bpi_{n+i \alpha}^{\Aone}(S^n \wedge \Gm^{\wedge m}) = \KMW_{m - i}
\end{equation}
provided $n \ge 2$ and $m \ge 1$. These calculations appear as \cite[Remark 6.42]{MorelmathbbAalgebraictopology2012}.

The following result of Morel is known as the ``unstable connectivity theorem''.
\begin{theorem}
 Let $n > 0$ be an integer and let $X$ be a pointed $(n-1)$- connected object in $\Spaces$. Then its
 $\Aone$-localization is simplicially $(n-1)$-connected.
\end{theorem}

\section{Application to vector bundles} \label{sec:VB}

\subsection{$\Aone$-classifying spaces}

A recurring distinction in $\Aone$-homotopy theory is the existence of two different notions of classifying space of an
algebraic group scheme $G$. In both cases, one wishes to begin with a contractible object $E G$ on which $G$ acts,
and then to form the quotient $(E G)/G$. The difference arises because there are two distinct notions of quotient.

The \textit{simplicial} or \textit{Nisnevich} classifying space $B G$ is constructed by taking a contractible object
$EG$ on which $G$ acts freely, and then forming the quotient $BG=(EG)/G$ in the category of simplicial Nisnevich
sheaves. Standard homotopical devices for the construction of $BG$, such as found in
\cite{MayClassifyingspacesfibrations1975} for instance, will invariably produce this classifying space, up to homotopy.

One may produce a new space, $B_{\et} G$, by taking a (derived)
pushforward of the pullback of $BG$ along the morphism of sites
\[ (\Sm_k)_{\et} \to (\Sm_k)_{\Nis}, \] see \cite[Section 4.1]{morelMathbbAHomotopyTheory1999}. 

The \textit{geometric} classifying space, $B_{gm}G$, is constructed in \cite[Section
4.2]{morelMathbbAHomotopyTheory1999}, using what they term an \textit{admissible gadget}, which is a generalization of
the construction of \cite{TotaroChowringclassifying1999}. The construction of \cite{morelMathbbAHomotopyTheory1999}
applies to \'etale sheaves of groups over a base $S$, but for the sake of the exposition here, we restrict to a reductive
algebraic group $G$ over a field $k$. While $EG$ cannot be constructed as a variety, one may construct a sequence
$U_i \hookrightarrow U_{i+1}$ of increasingly highly-$\Aone$-connected varieties on which $G$ acts freely and
compatibly. The construction is accomplished by considering a suitable family of larger and larger representations of
$G$ on affine spaces over $k$, and discarding the locus where the action of $G$ is not free. The quotients $U_i/G$ may
then be taken in the category of algebraic spaces, or if one is particularly careful, in schemes, and $B_{gm} G$ is then
defined to be the colimit of the ind-algebraic-space $\{U_i/G\}$. The quotient that is constructed here is `geometric', in
that it agrees with classically-existing notions of quotient of one variety by another, as in e.g.,
\cite{MumfordGeometricinvarianttheory1994}; it is also the quotient in the big site of \'etale sheaves. It is then the
case, \cite[Proposition 4.2.6]{morelMathbbAHomotopyTheory1999}, that $B_{\et} G$ is $\Aone$-equivalent to $B_{gm} G$.

If $G$ is a sheaf of \'etale group schemes, e.g., if $G$ is representable, and if $\Hoh^1_{\Nis}(\cdot, G) \iso
\Hoh^1_{\et} (\cdot , G)$, then the constructions of $B_{\et}G$ and $BG$ above are naturally simplicially equivalent,
\cite[Lemma 4.1.18]{morelMathbbAHomotopyTheory1999}. In particular, in the case of the \textit{special} algebraic groups
of \cite[Section 4.1]{SerreEspacesfibresalgebriques1995}, including $\GL_n$, $\SL_n$ and $\Sp_n$, all notions of
classifying space considered above are $\Aone$-weakly equivalent.

In contrast, in the case of a nontrivial finite group, both $BG$ and $B_{\et} G$ are $\Aone$-local objects already,
\cite[Section 4.3]{morelMathbbAHomotopyTheory1999}, and one may easily construct a $G$-torsor $\pi: U \to U/G$ where
both $U$ and $U/G$ are smooth $k$ schemes but where $\pi$ is not Nisnevich-locally trivial, it follows from \cite[Lemma
4.1.8]{morelMathbbAHomotopyTheory1999} again that for a nontrivial finite group, $BG\to B_{\et} G$ is never an
$\Aone$-equivalence.

In the study of vector bundles, or equivalently of $\GL_n$-bundles, the spaces $B_{\gm} \GL_n$ appearing above are
approximated by the ordinary Grassmanians, $\Gr_n(\A^{r})$. The notation $\Gr_n$ is adopted for the colimit as $r \to
\infty$, and the $\Aone$-weak equivalences
\[ \Gr_n \weq_{\Aone} B \GL_n \]
can be interpreted as a relation between a geometric construction on the left and a homotopical construction on the right.

\subsection{Classification of Vector Bundles}

The following metatheorem lies at the heart of the applications of $\Aone$-homotopy to the classification of vector bundles:
\begin{theorem}
 Let $k$ be a sufficiently pleasant base ring, and let $r \in \mathbb{N}$. The functor assigning to a smooth affine
 $k$-scheme $X$ its set of
 rank-$r$ vector bundles, denoted $X \mapsto V_r(X)$, is represented in the $\Aone$-homotopy category by the infinite
 Grassmannian of $r$-planes, $\Gr_n$. That is, there is a natural bijection of sets
 \begin{equation}
 \label{eq:1}
 V_r(X) \leftrightarrow \Mor_{\sh H(k)}(X, \Gr_r).
 \end{equation}
\end{theorem}

This metatheorem was proved by Morel as \cite[Theorem 8.1]{MorelmathbbAalgebraictopology2012} in the case where
$n \neq 2$and subject to
the running assumptions of \cite{MorelmathbbAalgebraictopology2012}, namely, that
$k$ be an infinite perfect field.\footnote{Although the case
$n=2$ can presumably be proved using similar methods, the work on the subject has not been published.} Subsequently this was vastly generalized, to the case where
$k$ is a regular noetherian ring over a Dedekind domain with perfect residue fields---this includes, in particular, all
fields---and allowing the case
$n=2$ in \cite[Theorem 5.2.3]{asokAffineRepresentabilityResults2017}.

We remark in passing that the restriction to smooth affine $k$-schemes cannot be weakened. Once the ``affine'' hypothesis is dropped, the set of isomorphism classes of vector bundles
is no longer $\Aone$-invariant, and so cannot be represented in the $\Aone$-homotopy category. One can consult
\cite{AsokVectorbundlescontractible2008} for a proof that $V_r(\PP^1) \to V_r(\PP^1 \times \Aone)$ is not an isomorphism
or that there exists a smooth quasiaffine $X \weq_{\Aone} \pt$ such that $V_r(X) \neq V_r(\pt)$.

One can use the metatheorem to deduce facts about $V_r(X)$. The main computational tool is the existence of a
Postnikov tower in $\Aone$-homotopy theory; which reduces the calculation of maps $X \to \Gr_n$, where $X$ is a smooth
affine variety, to a succession of lifting problems. It should be noted that $\bpi_1^{\Aone}(\Gr_n) \iso \Gm$, and this
$\bpi_1$ acts nontrivially on the higher homotopy sheaves of $\Gr_r$, so that the lifting problems one encounters for
$\GL_n$ are $\Gm$-twisted lifting problems. The essentials of the twisted sheaf-theoretic obstruction are set out in \cite[Section
6]{Asokcohomologicalclassificationvector2014} and \cite[Appendix B]{MorelmathbbAalgebraictopology2012}.

Once the obstruction theory has been set up, two related problems remain: the calculation of
$\pi_i^{\Aone}(\Gr_n)$ for various values of $i$---including the $\Gm$-action---and the interpretation of the resulting
lifting problems.

For the calculational problem, the following is known
\[ \xymatrix{ \A^n - \{0 \} \ar[r] & \Gr_{n-1} \ar[r] & \Gr_n } \]
is an $\Aone$-homotopy fiber sequence, \cite[Remark 8.15]{MorelmathbbAalgebraictopology2012}. Let $X=\Spec A$ denote a smooth
affine $k$-scheme and let $V$ denote a rank $n$ vector bundle on $X$, i.e., a projective module. The ``splitting
problem'' is to determine necessary and sufficient conditions for there to be an isomorphism $V \iso V' \oplus A$.

\subsection{The splitting problem}\label{Section:splittingproblem}

Since $\A^n - \{0\}$ is $n-2$-connected, by results of \cite[Chapter 6]{MorelmathbbAalgebraictopology2012}, it follows
by induction that if the dimension of $d< n$, then one can always split $V \iso V' \oplus A$. This gives a ``geometric''
argument for a result of Serre, \cite{SerreModulesprojectifsespaces1958}. 

More interesting is what happens when the obstruction-theory problem is not trivial for dimensional reasons. The problem
is slightly easier if all vector bundles appearing are assumed to be equipped with a trivialization of their determinant
bundles, so that $B\SL_n$ may be substituted for $\Gr_n = B \GL_n$. The first obstruction to the splitting problem is an
\textit{Euler class} class $\tilde c_n (V) \in \Hoh^n_{\Nis}(X; \KMW_n)$, as calculated in \cite[Theorem
8.14]{MorelmathbbAalgebraictopology2012}, arising directly from the obstruction theory in $\Aone$-homotopy theory, and
if the dimension of the base, $d$, is equal to the rank $n$ of the bundle, then this is the only obstruction.

Related to the above, Nori gave a definition of an ``Euler class group'' $E(A)$ of a noetherian ring $A$,
appearing in \cite[Section 1]{MandalEulerclassescomplete1996a} and of the Euler class $e(P)$ of a projective module, and
Bhatwadekar \& Sridharan \cite{BhatwadekarEulerClassGroup2000} established that the vanishing of the Euler class is
precisely equivalent to the splitting off of a trivial summand of $P$. There exists a surjective map $\eta: E(A) \to
\widetilde{CH}^n(\Spec A)$, where $A$ is a smooth $k$-algebra of dimension $n$, \cite[Chapitre
17]{faselGroupesChowWitt2008}, taking one Euler class to the other. According to recent, not yet published, work of Asok
and Fasel \cite{AsokEulerclassgroups2016}, the map $\eta$ is an isomorphism for smooth affine varieties over an infinite
perfect field of characteristic different from $2$.

\subsection{Other results on vector bundles}

Obstruction theory in $\Aone$-homotopy theory, applied to the representing space for $V_r$, as established in \eqref{eq:1},
is a powerful technique for deducing results about vector bundles.

The main result of \cite{Asokcohomologicalclassificationvector2014} is the following:
\begin{theorem} \label{th:SmAff3}
 Suppose $X$ is a smooth affine 3-fold over an
algebraically closed field $k$ having characteristic unequal to 2. The map assigning to a rank-2
vector bundle $E$ on $X$ the pair $(c_1(E), c_2(E))$ of Chern classes gives a pointed bijection 
\[ V_2(X) \overset{\iso}{\to} \Pic(X) \times CH^2 (X).\]
\end{theorem}
The surjectivity of the map above had previously been established in \cite{KumaAlgebraicCyclesVector1982}. The method of
proof is a calculation---at least in part---of the first three nontrivial $\Aone$-homotopy groups of $\GL_n$. Of these,
two are calculated by Morel in \cite{MorelmathbbAalgebraictopology2012}. Specifically $\bpi^{\Aone}_0(\GL_n) \iso \Gm$,
for all $n \ge 0$: the reason being that $\SL_n$ is $\Aone$-connected \cite[6.52]{MorelmathbbAalgebraictopology2012}, so
that
\[ 1 \to \SL_n \to \GL_n \to \Gm \to 1 \]
is an $\Aone$-homotopy fibre sequence. It is moreover split, and $\Gm$ is strongly $\Aone$-invariant, so it follows that for $n \ge 1$
\[ \bpi^{\Aone}_i(\GL_n) = \begin{cases} \Gm \text{\quad if $i=0$,} \\ \bpi^{\Aone}_i(\SL_n) \text{\quad
 otherwise.} \end{cases} \]
It is then also possible to calculate $\bpi^{\Aone}_1(\GL_2) = \bpi^{\Aone}_1(\SL_2) = \bpi^{\Aone}_1(\A^2-\{0\}) \iso
\KMW_2$---this belonging in the family of calculations of $\bpi^{\Aone}_n(S^{n+m\alpha})$ of
\cite{MorelmathbbAalgebraictopology2012}---the reason the groups $\bpi_1^{\Aone}$ are sometimes exceptional is that they are
not known in general to be abelian, but here one is calculating a higher homotopy sheaf of a group, so the
Eckmann--Hilton argument applies and $\bpi_1^{\Aone}(\SL_2)$ is abelian.

The sheaf $\KM_2$ relates to classical invariants of vector bundles via the ``formula of Rost'': 
\[ \Hoh_{\Nis}^n(X, \KM_n) = \Hoh_{\Zar}^n(X, \KM_n) = CH^n(X), \]
the first identity being \cite{MorelmathbbAalgebraictopology2012} and the second being \cite[Corollary
6.5]{RostChowgroupscoefficients1996}. For a quadratically closed field, the result also holds for $\KMW_n$, see
\cite{Asokcohomologicalclassificationvector2014}.

Knowing the first two nonvanishing homotopy sheaves of $B\GL_2$, along with the action of $\bpi_1^{\Aone}(B\GL_2) = \Gm$
on $\bpi_2^{\Aone}$, allows one to construct the second stage of the Postnikov tower, $B\GL_2^{(2)}$. The action of
$\Gm$ does have to be taken into account, so the calculation is more involved than simply the observations
\[ [X, B \Gm]_{\Aone} = \Pic(X), \quad [X, K(\KMW_2(\sh L), 2)]_{\Aone} = CH^2(X), \]
for a smooth $X$---here $\KMW_2(\sh L)$ is a twisted form of $\KMW_2$---, but over a quadratically closed field the
ultimate calculation is the same. For a $3$-fold, there is no further obstruction to lifting to a map
$[X, B \GL_2]_{\Aone}$, but what there might be \textit{a priori} is a choice of different liftings. The set of such
choices is $\Hoh^3(X, \bpi^{\Aone}_3(B\GL_2)(\sh L))$, and \cite{Asokcohomologicalclassificationvector2014} shows that
this is trivial over an algebraically closed field. It is worth remarking that the sheaf $\bpi^{\Aone}_3(B\GL_2)$ in
question---at least in the untwisted form---is isomorphic to $\bpi^{\Aone}_2(S^{1+2\alpha})$, that is, it corresponds to
the stable $1$-stem.

A number of results in a similar vein to the above can be established by the methods outlined. It is a corollary of Theorem
\ref{th:SmAff3} that a rank-2 vector bundle on a smooth affine 3-fold over an algebraically closed field has a trivial
line-bundle summand if and only if $c_2$ vanishes. In \cite{AsokSplittingvectorbundles2015}, the generalization to
dimension $n$ is named ``Murthy's Splitting Conjecture'': 

{\small ``On an affine variety of dimension $n$ over an algebraically
closed field, a rank $n-1$ vector bundle admits a trivial summand if and only if $c_{n-1}$ vanishes.''}

The conjecture is proved for $n=4$, again by analysis of a group in the $1$-stem: $\bpi_3^{\Aone}(\A^3 - \{0\})$. In
spite of considerable attention, \cite{AsokmetastablerangemathbbA2013}, \cite{AsokexplicitmathrmKOdegree2017},
\cite{Asoksimplicialsuspensionsequence2017}, the groups $\bpi_n^{\Aone}(\A^n- \{0\})$ have not been calculated for
$n \ge 4$, and the problem remains open.

The groups $\GL_n$ and $\SL_n$ are the most studied groups, corresponding to the study of vector bundles, either
\textit{per se} or with an orientation. Other groups have also been studied. For instance, Wendt shows in
\cite{WendtRationallytrivialtorsors2011} that if $G$ is a smooth split reductive group over an infinite field $k$, and if $E
\to B$ is a $G$-torsor satisfying certain local triviality conditions---local triviality in the Zariski topology and
triviality over the basepoint of $B$---then the fibre sequence
\[ G \to E \to B \]
is an $\Aone$-homotopy fibre sequence. In \cite{AsokAffinerepresentabilityresults2018}, the authors establish, again
over an infinite field, that for an isotropic reductive $k$ group $G$, the functor $\Hoh_{\Nis}^1(\cdot, G)$ of
isomorphism classes of Nisnevich-trivial $G$-torsors is represented in the $\Aone$-homotopy category by $BG$. These
results generalize to a wide class of groups the results of Morel's that are foundational to the study of vector
bundles. Unfortunately for the applicability of the theory, however, for many algebraic groups $G$ of interest, the set of
torsors of interest is $\Hoh_{\et}^1(\cdot, G)$ rather than $\Hoh_{\Zar}^1(\cdot, G)$ or $\Hoh_{\Nis}^1(\cdot, G)$; therefore the theory is most
interesting in the case of special algebraic groups in the sense of \cite{SerreEspacesfibresalgebriques1995}---these
groups include $\SL_n$, $\GL_n$ and $\Sp_n$. In \cite{AsokGenericallysplitoctonion2017}, the same authors apply their techniques to
$G_2$, where $\Hoh^1_{\Nis}(\cdot, G_2)$ classifies ``split octonion algebras''.

The means by which Asok--Hoyois--Wendt establish their representability result is a strong general representability result: working
over a general quasicompact, quasiseparated base scheme, and letting $\mathcal F$ denote a presheaf such that $\pi_0(\sh F)$ is an
$\Aone$-invariant functor and such that $\sh F$ satisfies a certain ``affine Nisnevich excision'' property---a Mayer--Vietoris style
property---then $\pi_0(\sh F)(U) = [U , \sh F]_{\Aone}$ for affine $S$-schemes, \cite[Theorem 5.1.3]{asokAffineRepresentabilityResults2017}.
If one restricts to an infinite base field $k$, then \cite{AsokAffinerepresentabilityresults2018} shows that if $G$ is an isotropic
reductive group scheme, then the simplicial classifying space $B G$---or, more correctly, a Nisnevich fibrant replacement---satisfies these
hypotheses. 

\subsection{Naive $\Aone$-homotopy}

If $\sh F : \cat{Ring} \to \cat{Set}$ is a functor defined on (affine) rings, or some subcategory of rings, then it is possible to
define \textit{na\"ive $\Aone$-homotopy classes} of elements in $F(R)$. One considers the two evaluation maps $R[t] \to R$ given by
evaluation at $t=0$ and $t=1$. A na\"ive homotopy between $a$ and $b$ in $F(R)$ is an element $H \in F(R[t])$ such that under the two
evaluation maps $H \mapsto a$ and $H \mapsto b$. Of course, there is no reason in general to this should be an equivalence relation,
so one takes the equivalence relation generated by it.

The definition above also applies to contravariant functors on schemes, but it is on affine schemes that it works best. When $\sh F$ is
representable, say by a scheme $X$, then na\"ive homotopy on $\sh F( R)$ in the sense above recovers the theory of equivalence classes of
maps $\Spec R \to X$ under the relation generated by declaring the two restrictions of a map $H : \Spec R[t] \to X$ at $t=0$ and $t=1$ to be
equivalent.

It is worth noting that the na\"ive $\Aone$-equivalence classes of maps $\Spec R \to X$ are exactly the path components of the $\Aone$
singular construction $\Sing^{\Aone} X$ evaluated at $R$, i.e., $\pi_0( \Sing^{\Aone} X(R))$. In general, $\Sing^{\Aone} X$ is not an
$\Aone$-local object, so that na\"ive $\Aone$-homotopy is not the same as genuine $\Aone$-homotopy. This is perhaps the great challenge in unstable
$\Aone$-homotopy---the na\"ive, geometrically appealing construction does not coincide with the genuine construction possessed of good
theoretical properties. One consequence of the work of \cite{asokAffineRepresentabilityResults2017},
\cite{AsokAffinerepresentabilityresults2018}, however, is that in a number of cases, genuine and na\"ive
$\Aone$-homotopy agree.

For an example where the na\"ive and genuine $\Aone$-homotopy classes of maps differ, we refer to
\cite{cazanaveHomotopyClassesRational2008}. In this paper, the na\"ive $\Aone$-homotopy classes of maps
$\PP^1 \to \PP^1$ are calculated. In order to put this in context, we should mention that in
\cite[Theorem 7.36 and Remark 7.37]{MorelmathbbAalgebraictopology2012}, Morel calculated
\begin{equation}\label{unstableA1homotopyclassesP1toP1} [\PP^1 , \PP^1]_{\Aone} = \GW(k) \times_{k^\times/{k^\times}^2} k^\times. \end{equation}
It is to be remarked that although $\PP^1$ is a motivic sphere, this calculation lies outside the main body of
calculations of self-maps of spheres, since $\PP^1 \weq S^1 \wedge \Gm$ is not $\Aone$-simply-connected. The main
theorem of \cite{cazanaveHomotopyClassesRational2008} is that the na\"ive homotopy classes of maps of spheres make up
the submonoid of $\GW(k) \times_{k^\times/{k^\times}^2} k^\times$ corresponding to elements in $\GW(k)$ of positive
degree, so that the map from the set of na\"ive $\Aone$-homotopy classes of maps $\PP^1 \to \PP^1$ to the set of all
$\Aone$-homotopy classes is, in fact, a group completion.

The following is a part of \cite[Theorem 2.3.2]{AsokAffinerepresentabilityresults2018}. As before, the base $S$ is quasicompact and
quasiseparated.
\begin{theorem}
 Suppose $G$ is a finitely presented smooth $S$-group scheme and that $H$ is a finitely presented closed subgroup scheme such that $G/H$ is
 an $S$-scheme, and such that $G \to G/H$ is Nisnevich locally split, and suppose that $\Hoh^1_{\Nis}(\cdot, G)$ and $\Hoh^1_{\Nis}(\cdot,
 H)$ are $\Aone$-invariant on affine schemes, then for an affine scheme $\Spec R$, the natural map
 \[ \pi_0(\Spec^{\Aone}(G/H)(R)) \to [\Spec R, G/H]_{\Aone} \]
 is a bijection.
\end{theorem}
That is, for $G/H$, provided the source is affine, the na\"ive and genuine $\Aone$-homotopy classes of maps agree. The restrictions on $G$,
$H$ are not perhaps as inconvenient as might be feared, since, in contrast to the case of the \'etale topology, the Nisnevich cohomology $\Hoh^1_{\Nis}(\cdot, G)$ is
liable to be $\Aone$-invariant. In fact, in \cite[Theorem 3.3.6]{AsokAffinerepresentabilityresults2018}, the authors establish $\Aone$
invariance for all isotropic reductive group schemes over infinite fields. The more severe restriction in the theorem is that $G \to
G/H$ should be Nisnevich locally split, but nonetheless, the theorem provides a wide range of notable spaces for which na\"ive and genuine
$\Aone$-homotopy are the same.

 \bibliographystyle{alpha} \bibliography{BenW_Standard_BibTeX,Chapter}

\end{document}